\documentstyle[11pt,righttag]{amsart}

\input{psfig}

\catcode`\@=11

\long\def\@savemarbox#1#2{\global\setbox#1\vtop{\hsize\marginparwidth 
  \@parboxrestore\tiny\raggedright #2}}
\newcommand\lref[1]{\ref{#1}%
\@ifundefined{r@DisplaY #1}{}{(#1)}}

\newcommand\fakelabel[2]{\@bsphack\if@filesw {\let\thepage\relax
   \newcommand\protect{\noexpand\noexpand\noexpand}%
\xdef\@gtempa{\write\@auxout{\string
      \newlabel{#1}{{#2}{\thepage}}}}}\@gtempa
   \if@nobreak \ifvmode\nobreak\fi\fi\fi\@esphack}

\catcode`\@=12

\def\Empty{}
\newcommand\oplabel[1]{
  \def\OpArg{#1} \ifx \OpArg\Empty {} \else
  	\label{#1}
  \fi}
\newtheorem{theoremSt}{Theorem}[section]

\newtheorem{exampleSt}[theoremSt]{Example}
\newtheorem{exerciseSt}[theoremSt]{Exercise}

\newcommand\MakeStEnv[1]{
  \newenvironment{#1}[1]{
  \begin{#1St} \oplabel{##1}%
  \global\def\CrntSt{\thetheoremSt}%
}{ 
  \end{#1St} }
  \newenvironment{#1+}[1]{
  \begin{#1St} \label{##1}%
  \label{DisplaY ##1}%
  \global\def\CrntSt{\thetheoremSt}%
  \def\Labl{##1}\ifx\Labl\Empty{} \else {\em (\Labl)\,}\fi%
}{ 
  \end{#1St} }
}
\MakeStEnv{theorem}
\MakeStEnv{corollary}
\MakeStEnv{proposition}
\MakeStEnv{lemma}
\MakeStEnv{define}
\MakeStEnv{conjecture}

\newcommand\restate[3]{
\medskip\par\noindent
{\bf #1 \ref{#2}} 
{\it #3}
\par\medskip
}

\long\def\realfig#1#2#3{
\begin{figure}[htbp]
\centerline{\psfig{figure=#2}}
\caption[#1]{#3}
\oplabel{#1}
\end{figure}}
			 		
\newlength{\saveu}

\newenvironment{proof}[1]{%
  \def\PfArg{#1}
  \ifx\PfArg\Empty
  	\edef\PfArg{\CrntSt}  \fi
 {\parindent=0pt\startproof{\PfArg}}%
}{ 
  \finishproof{\PfArg}
}
\newcommand{\startproof}[1]{%
\medbreak\mbox{}\noindent{\it Proof of #1:}%
}
\newcommand{\finishproof}[1]{ 
  \def\FPArg{#1}
  \ifx\FPArg\Empty
  	\newcommand\FPArg{\CrntSt}  \fi
  \smallbreak\noindent\makebox[\textwidth]{\hfill\fbox{\FPArg}}
  \medbreak\noindent
}

\newcommand{\bfheading}[1]{\par\smallskip\noindent {\bf #1}}

\newcommand\AAA{{\cal A}}

\newcommand\CC{{\cal C}}

\newcommand\MM{{\cal M}}
\newcommand\NN{{\cal N}}

\newcommand\PP{{\cal P}}
\newcommand\QQ{{\cal Q}}
\newcommand\RR{{\cal R}}

\newcommand\TT{{\cal T}}

\newcommand\half{{\textstyle{1\over2}}}

\newcommand\Area{\operatorname{Area}}

\newcommand\ep{\epsilon}

\newcommand\union{\cup}
\newcommand\intersect{\cap}
\newcommand\bbR{{\mathord{\text{I\kern-2pt R}}}}        
\newcommand\bbH{{\mathord{\text{I\kern-2pt H}}}}        

\newcommand\C{{\bold C}}

\newcommand\Z{{\bold Z}}
\newcommand\R{{\bold R}}

\newcommand\Hyp{{\bold H}}

\newcommand\SL[1]{\text{SL}_{#1}}

\newcommand\bigrightarrow[1]{\hbox to #1{\rightarrowfill}}
\newcommand\bigleftarrow[1]{\hbox to #1{\leftarrowfill}}

\newcommand\boundary{\partial}
\newcommand\semidir{\mathrel{\hbox{\vrule depth-.03ex height1.1ex\kern-0.15em$\times$}}}

\newcommand\til{\widetilde}

\newcommand\gesim{\succ}
\newcommand\lesim{\prec}
\newcommand\simle{\lesim}
\newcommand\simge{\gesim}
\newcommand{\simmult}{\asymp}

\numberwithin{equation}{section}

\catcode`\@=11

\def\@sect#1#2#3#4#5#6[#7]#8{%
\ifnum #2>\c@secnumdepth
   \let\@svsec\@empty
 \else
   \refstepcounter{#1}%
\edef\@svsec{\ifnum#2<\@m
             \@ifundefined{#1name}{}{\csname #1name\endcsname\ }\fi
\noexpand\rom{\S \csname the#1\endcsname.}\enspace}\fi
 \@tempskipa #5\relax
 \ifdim \@tempskipa>\z@ 
   \begingroup #6\relax
   \@hangfrom{\hskip #3\relax\@svsec}{\interlinepenalty\@M #8\par}%
   \endgroup
      \iffalse\else\csname #1mark\endcsname{%
     \ifnum \c@secnumdepth >#2\relax\@svsec\fi#7}\fi
\ifnum#2>\@m \else
       \addcontentsline{toc}{#1}%
{\ifnum #2>\c@secnumdepth \else
             \protect\numberline{%
               \ifnum#2<\@m
               \@ifundefined{#1name}{}{\csname #1name\endcsname\ }\fi
               \csname the#1\endcsname.}\fi
           #8}%
     \fi
 \else
  \def\@svsechd{#6\hskip #3\@svsec
    \@ifnotempty{#8}{\ignorespaces#8\unskip
       \ifnum\spacefactor<1001.\fi}%
        \ifnum#2>\@m \else
          \addcontentsline{toc}{#1}%
            {\ifnum #2>\c@secnumdepth \else
              \protect\numberline{%
                \ifnum#2<\@m
                \@ifundefined{#1name}{}{\csname #1name\endcsname\ }\fi
                \csname the#1\endcsname.}\fi
             #8}\fi}%
 \fi
\@xsect{#5}}
\def\@evenhead{\defaultfont\small\sc
      \rlap{\thepage}\hfil 
      \expandafter{\sh@rttitle}\hfil}%
\def\@oddhead{\defaultfont\small \sc \mbox{}\hfil
      \expandafter{\rightmark}\hfil 
      \llap{\thepage}}%

\def\section{\@startsection{\@string\section}%
        1\z@{9\p@\@plus12\p@}{6\p@}%
        {\centering\defaultfont\large\bf}}

\catcode`\@=12

\begin{document}

\title[Product regions in Teichm\"uller space]{Extremal length
estimates and product regions in Teichm\"uller space}
\author{Yair N. Minsky}
\address{Mathematics Department and IMS, SUNY Stony Brook, NY 11794}
\email{yair@math.sunysb.edu}

\date{September 1, 1994}
\thanks{This work was partially supported by an NSF postdoctoral
fellowship.} 

\maketitle

\begin{abstract}
We study the Teichm\"uller metric on the Teichm\"uller space of a surface 
of finite type, in regions where the injectivity radius of the surface is 
small. The main result is that in such regions the Teichm\"uller metric 
is approximated up to bounded additive distortion by 
the sup metric on a product of lower dimensional spaces. 
The main technical tool in the proof is the 
use of estimates of extremal lengths of curves in a surface based on the 
geometry of their hyperbolic geodesic representatives.
\end{abstract}

\section{Introduction}
There is a longstanding but imperfect analogy between the geometry of the
Teichm\"uller space $\TT(S)$ of a surface and that of a complete, negatively
curved space. For example, $\TT(S)$ admits a boundary at infinity
similar to that of hyperbolic space (See e.g.
\cite{travaux,masur:boundaries}), and the Teichm\"uller geodesic flow on its
quotient, the moduli space, is ergodic (see \cite{masur:geoflow}).
This paper studies one of the strong ways in which this analogy
fails, namely the existence of large regions in the space which are
closely approximated by products of lower-dimensional metric spaces. 

There are several natural metrics on $\TT(S)$; we will consider
throughout only the {\em Teichm\"uller metric} (\S\ref{kerckhoff theorem}).
The first indication of positively curved behavior of the
Teichm\"uller metric came from Masur \cite{masur:teichgeo}, in which
examples were given of geodesic rays with a common basepoint, which
remain a bounded distance apart for all time. This contradicts
non-positive curvature on small scales. More recently, Masur and Wolf
\cite{masur-wolf} have given examples of geodesic triangles which
fail, arbitrarily badly, the ``thin triangle'' condition for
hyperbolicity in the sense of Gromov and Cannon
\cite{cannon:negative,c-d-p,gromov:hypgroups}; this implies that the 
large-scale behavior of the metric is not negatively curved.  The main
result of this paper describes regions of $\TT(S)$ where the large
scale behavior exhibits some characteristics of positive
curvature.

The main theorem can be summarized as follows (a complete statement
appears in \S \ref{product structure}). Let $\TT(S)$ be the
Teichm\"uller space of a surface $S$ of finite 
type, endowed with the Teichm\"uller metric. Let $\gamma =
\gamma_1,\ldots,\gamma_k$ be a system of disjoint, 
homotopically distinct simple closed curves on $S$, and let
$Thin_\ep(S,\gamma)$ denote the set of $\sigma\in\TT(S)$ for which
$\ell_\sigma([\gamma_i])\le \ep $ for all $i$, where
$\ell_\sigma([\gamma_i])$ denotes 
hyperbolic length of the $\sigma$-geodesic representative of $\gamma_i$.
Let $X_\gamma$ denote the product space 
$\TT(S\setminus\gamma) \times \Hyp_1\times\cdots\times \Hyp_k$, where
$S\setminus\gamma$ is 
considered as a punctured surface, and each $\Hyp_i$ is a copy of the
hyperbolic plane. Endow $X_\gamma$ with the {\em sup metric}, $d_X =
\max\{d_{\TT(S\setminus\gamma)},d_{\Hyp_1},\ldots,d_{\Hyp_k}\}$, of the 
metrics on the factors. 
 
\restate{Theorem}{Product structure}{
The Fenchel-Nielsen coordinates on $\TT(S)$ give rise to 
a natural homeomorphism $\Pi:\TT(S)\to X_\gamma$, 
and for $\ep$ sufficiently small this homeomorphism restricted to
$Thin_\ep(S,\gamma)$ distorts distances by a bounded {\em additive
amount}. That is, 
$$
|d_{\TT(S)}(\sigma,\tau) - d_X(\Pi(\sigma),\Pi(\tau))| \le c
$$
For $\sigma,\tau\in Thin_\ep(S,\gamma)$, and $c=c(\chi(S),\ep)$.
}

Note in particular that the distortion of $\Pi$  is additive and not
multiplicative, as is more commonly the case with quasi-isometries. For
this reason it is significant that we endow $X_\gamma$ with the sup
metric, and not the $L_2$ metric; these are equivalent up to bounded
multiplicative, but not additive, distortion. 
A product of spaces with the sup metric exhibits not just
non-negative, but actually positive curvature characteristics in the
large. See \S\ref{remarks} for a fuller discussion of this point. 

\medskip 

The idea of decomposing Teichm\"uller space as a product using length
and twist parameters of simple closed curves dates back to
Fenchel-Nielsen \cite{fenchel-nielsen} (see Abikoff \cite[Chap.
II]{abikoff} and Fricke-Klein \cite{fricke-klein}). Such 
decompositions have 
particular significance when the curves in question are short, and
this has been used by several authors to study the degeneration of
Riemann surfaces towards ``noded'' surfaces, in which a short curve
has been replaced by two punctures. From a complex-analytic and
algebraic-geometric point of view, such degenerations 
have been investigated by e.g.
Bers \cite{bers:degenerating},
Earle-Marden \cite{earle-marden}, Kra \cite{kra:horocyclic}, 
Maskit \cite{maskit:decomposition,maskit:moduli},
Masur \cite{masur:wpmetric} and 
Wolpert \cite{wolpert:cutpaste,wolpert:plumbing,wolpert:wpgeom}.
In all these cases, the emphasis is on careful local analysis,
concentrating on the complex structure of the moduli space in a
neighborhood of a noded surface, or on 
differential geometric structure such as the Weil-Petersson metric
in the case of \cite{masur:wpmetric,wolpert:wpgeom}, or a precise
expression for the hyperbolic metric 
on the surface in \cite{wolpert:plumbing}.
By contrast, our approach ignores most local analytical details
and concentrates on large-scale estimates in the Teichm\"uller space.

\bfheading{Extremal length estimates.}
The estimates of Teichm\"uller distance necessary for theorem
\ref{Product structure} are obtained using Kerckhoff's theorem (see
\S\ref{kerckhoff theorem}), which relates Teichm\"uller distance  to
ratios of extremal lengths. Therefore our main tool is a theorem that
allows us to estimate extremal lengths of arbitrary curves on a
surface, in terms of the behavior of their hyperbolic geodesic
representatives. 

The relationship between hyperbolic length $\ell_\sigma$ and extremal
length $\lambda_\sigma$
(see \S\ref{length defs}) of curves on a Riemann surface $(S,\sigma)$ is not
completely straightforward. For very short curves, $\ell_\sigma$ and
$\lambda_\sigma $ are nearly proportional.
For long curves, the 
relation depends on the geometry of the surface. In particular
an upper  bound for $\lambda_\sigma$ is given by $\lambda_\sigma \le
(1/2)\ell_\sigma e^{\ell_\sigma/2}$ (see 
Maskit \cite{maskit:comparison}), and this is close to sharp in general. But
on the other hand $\lambda_\sigma \ge \ell_\sigma^2/2\pi|\chi(S)|$ is
the general lower bound that comes from the definition of extremal
length, and for any {\em fixed} conformal structure $\sigma$  one can show (see
e.g. \cite{minsky:slowmaps} and lemma \ref{extremal is hyperbolic for thick})  
that it is the lower bound which estimates $\lambda_\sigma$
to within a bounded factor.

Theorem \ref{main estimate} gives an estimate for $\lambda_\sigma$
which is correct, up to a bounded factor that depends only on the
topological type of the surface, for all cases. The main observation
is that, after decomposing a Riemann surface $(S,\sigma)$ into thick
and thin parts, or what we call here an $(\ep_0,\ep_1)$ collar
decomposition (see \S\ref{thick-thin}), each component of the 
decomposition contributes individually to the extremal length of a
curve $\beta$ in $S$. The contribution $\lambda_{P,\sigma}(\beta)$ of
a component $P$ of the thick part (see definition (\ref{P length def})
and lemma \ref{extremal is hyperbolic for thick})
is approximately $\ell_\sigma^2(\beta^\sigma\intersect P)$, where
$\beta^\sigma$ is the $\sigma$-geodesic representing $\beta$. A
component $A$ of the thin part, which is an annulus, contributes a quantity
$\lambda_{A,\sigma}(\beta)$ which depends on the modulus of $A$ and
the amount of ``twisting'' that $\beta$ does in $A$ (see
\S\ref{twisting section} and definition (\ref{A length def})).
The theorem then states: 

\restate{Theorem}{main estimate}{
Let $S$ be a surface of finite type with
a standard hyperbolic metric $\sigma$ in which
boundary lengths are at most $\ell_0$, and let $\QQ$ be the set of 
components of the $(\ep_0,\ep_1)$ collar decomposition.
Then, for any curve system $\beta$ in $S$, the extremal length 
$      \lambda_{S,\sigma}(\beta)$ is approximated by
$$
                          \max_{Q\in\QQ} \lambda_{Q,\sigma}(\beta)
$$
up to multiplicative factor depending only on 
$\ep_0,\ep_1,\ell_0$ and the topological type of $S$.
}
\noindent
In specific cases this theorem recovers the above exponential and
quadratic bounds from \cite{maskit:comparison} and \cite{minsky:slowmaps}.

\bfheading{Proofs.}
The proof of 
theorem \ref{main estimate} is essentially a long exercise in
elementary hyperbolic geometry. Given a homotopy class of curves
$\beta$ in $S$ we must arrange it in a sufficiently well-spaced way so
as to obtain an embedded annulus which gives us an upper bound for
extremal length (lower bounds are easier). The part of the argument
which is not straightforward 
has to do with finding an appropriate metric, conformally equivalent
to the hyperbolic metric, in which to make the spacing construction.
The scaling factors, which vary among components of the thick-thin
decomposition, must depend on the way in which $\beta$ is distributed
on $S$. The construction of this metric is given in lemma
\ref{Scaling function} and the discussion that follows it. 
(Note that the metric we obtain is 
in some way mimicking the actual extremal metric, which is given by a
quadratic differential. However our construction is explicitly
related to the hyperbolic characteristics of $\beta$, which makes it
more suited to our purposes.)

The proof of theorem \ref{Product structure} is an almost direct
application of theorem \ref{main estimate} and Kerckhoff's theorem
\ref{Kerckhoff distance thm}, except for one point. In the passage
from a component of $S\setminus \gamma$ with (short) boundary to a
punctured surface, we must at some point compare extremal length
ratios $\lambda_\sigma/\lambda_\tau$ for properly embedded arcs with
endpoints on the boundary to similar 
ratios for simple closed curves in the interior. That the suprema of
these ratios are within a bounded factor of each other is shown in
lemma \ref{closed curves give ratio}. The main trick in this lemma is
start with a (possibly non-closed) curve $\beta$ that maximizes the {\em
reciprocal} ratio $\lambda_\tau/\lambda_\sigma$, and then find a
closed curve which is sufficiently close to being ``orthogonal''
to $\beta$. 

A different approach to the proof of theorem \ref{Product structure}
than the one we have taken might involve estimating Teichm\"uller
distances directly by constructing appropriate quasiconformal maps.
One might hope that, since a short curve corresponds to a high-modulus
annulus in the surface, a quasi-conformal map on the complement of
this annulus might be ``spliced'' across the annulus in a way that
takes the appropriate modulus and twist parameters into account, and
thus the quasi-conformal dilatation of the spliced map would be
estimated by the dilatations of the component maps (and the modulus
and twisting information in the annulus).  It is therefore of interest
that, in our experience, this approach actually runs into quite
serious difficulties, since ``splicing'' quasi-conformal maps is, in
general, not very easily done without great loss of control over
dilatation. In particular, estimates of dilatation of such spliced
maps tend to be, at best, polynomially related to the dilatations of
the component maps, and not linearly as we would need for theorem
\ref{Product structure}. It is for this reason that we have adopted
the extremal length approach, which circumvents these problems. 

\bfheading{Summary of contents.}
In section \ref{defs} we set our notation and recall some well-known
facts about conformal structures on surfaces, curves and their
lengths, and Teichm\"uller distance. The only proof we give is an
elementary one for lemma \ref{Kerckhoff-Royden for torus}.

In section \ref{twisting section} we define {\em
twisting numbers}, which are a convenient way to measure how many
times one geodesic twists around another in a hyperbolic surface, and
in particular inside a Margulis annulus. We
prove a few elementary facts about them, and relate them to
Fenchel-Nielsen twist parameters in lemma \ref{compare twists}. 

Section \ref{basic extremal} gives some well-known estimates of
extremal lengths, in particular an estimate in terms of hyperbolic
length in a component of the thick part (lemma \ref{extremal is
hyperbolic for thick}). Lemma \ref{adjust on P}
describes briefly how to space a curve evenly in a component of the
thick part, and this will be used in \S\ref{main length estimate}. At
the end of the 
section we define the annulus contribution $\lambda_{A,\sigma}$ to
extremal length, in terms of twisting numbers. 

Sections \ref{main length estimate} and \ref{product structure}
give the proof of theorems
\ref{main estimate}, and \ref{Product structure}, respectively. 
Section \ref{remarks} gives a brief discussion of the geometric
consequences of theorem \ref{Product structure}.

\bfheading{Acknowledgements.}
The author is grateful to the Long Island Rail Road for its
generosity in providing space and time for contemplation.

\section{Definitions and background}
\label{defs}

\subsection{Metrics and conformal structures on surfaces.}
Let $S$ be a smooth orientable surface of finite topological type, by 
which we mean that $S$ is homeomorphic to a connected, compact surface (possibly 
with boundary) minus a finite set, which we call punctures or cusps.

By {\em conformal structure} on $S$ we shall mean a conformal structure 
on the interior, for which 
each puncture has a neighborhood conformally equivalent to a punctured 
disk, and each boundary component has a neighborhood conformally 
equivalent to an annulus of finite modulus.
Let $\TT(S)$ denote the Teichm\"uller space of marked conformal structures 
on $S$ considered up to conformal homeomorphisms homotopic to the
identity (see \cite{abikoff}).
If $S$ has no boundary (but possibly punctures) then $\TT(S)$ is the
space of {\em analytically finite} marked conformal structures.

A conformal structure also determines a class of metrics on the
interior, all those that can be described in a local conformal coordinate as
a positive function times the Euclidean metric. In this paper we will
use metrics which are smooth or mildly singular (non-differentiable on
some curves), and which have a completion to the boundary of $S$.  The
conformal class of a metric $\sigma$ will be denoted $[\sigma]$.

Call $S$ {\em hyperbolic} if its Euler 
characteristic is negative, which implies
that it admits a complete finite-volume hyperbolic metric
in which the boundary is geodesic. We will call such a hyperbolic
metric on $S$  {\em standard}.
We say a Euclidean metric on an annulus is
standard if the annulus can be obtained isometrically from a Euclidean
rectangle by identifying a pair of opposite edges.

\subsection{Curve systems.}

A {\em curve system} on $S$ is a disjoint union of homotopically 
non-trivial properly embedded simple arcs and closed curves (for an arc, 
non-trivial means not deformable into $\boundary S$ with endpoints fixed, 
and for a closed curve it means not deformable to a point or a puncture). 
For a hyperbolic $S$, let $\CC(S,\boundary S)$ denote the set of 
homotopy classes of curve systems under homotopies that keep the 
endpoints (if any) on the boundary. (We admit non-simple curves to
the homotopy classes for convenience.)
Let $\CC(S)\subset\CC(S,\boundary S)$ denote 
the subset of curve systems represented by unions of simple closed curves. 
Note that $\CC(S)$ includes curves that are homotopic to components of 
$\boundary S$. These are called {\em peripheral}, and to exclude them we 
define $\CC_0(S)$ to be the subset of $\CC(S)$ consisting of curve 
systems with no peripheral components. Also define (for brief use in
\S 3) $\CC'(S,\boundary S)$ to be the set of homotopy classes of curve
systems which are allowed to terminate in punctures as well as
boundary components.

For an annulus $A$, we make a different definition: Let $\CC(A,\boundary 
A)$ denote the  set of curve systems on $A$ with equivalence under 
homotopies that {\em  fix} the  endpoints.

We call a subsurface $M\subset S$ {\em hyperbolic} if no component of
$\boundary M$ is homotopically trivial or homotopic to another
component of $\boundary M$. Equivalently,
if $\sigma$ denotes any smooth non-positively curved metric on $S$ in which
the boundary is geodesic (or convex), then $M$ is isotopic to a
subsurface $M^\sigma$ whose boundary is geodesic in $\sigma$. 

Given a hyperbolic subsurface $M$, there is a natural map 
$$\RR=\RR_{M}:\CC(S,\boundary S) \to \CC(M,\boundary M),$$
defined as follows: Let $\sigma$ be a standard hyperbolic metric on $S$. 
If $\gamma$ is a curve system on $S$ let $\gamma^\sigma$ denote
its representative of shortest $\sigma$-length (see \S \ref{length 
defs}). Let $\RR([\gamma])$ be the class
represented by $\gamma^\sigma\intersect M^\sigma$. It is not hard to
see that $\RR$ does not depend on the choice of $\sigma$. 


\subsection{Lengths and geodesic representatives.} \label{length defs}
For any metric $\sigma$ on $S$ and any curve $\alpha$ we denote by
$\ell_\sigma(\alpha)$ the length of $\alpha$ as measured in $\sigma$.
If $[\gamma]\in\CC(S,\boundary S)$ then $\gamma^\sigma$ denotes a
representative of minimal length -- if $\sigma$ is a standard
hyperbolic metric then
this is unique, and orthogonal to $\boundary S$. For an annulus $A$
where the homotopies must fix endpoints, 
there is  also a minimal
representative $\gamma^\sigma$, but it may not be orthogonal to
$\boundary S$. In both cases we may denote
$\ell_\sigma(\gamma^\sigma)$ by $\ell_\sigma([\gamma])$. 
Note also that any $\gamma\in\CC'(S,\boundary S)$ has a unique
geodesic representative $\gamma^\sigma$, though its length may be
infinite. 

\bfheading{Extremal length.}
For a metric $\sigma$ and any class $[\gamma]$ of curves in $S$,
we can define the extremal length
$\lambda_\sigma([\gamma])$ as
\begin{equation}
\label{extremal length definition}
\lambda_\sigma([\gamma]) = \lambda_{S,\sigma}([\gamma]) = 
\sup_{\sigma'\in[\sigma]}
{\displaystyle\inf_{\gamma'\in[\gamma]} 
				\ell^2_{\sigma'}(\gamma') \over
				\Area(\sigma')}.
\end{equation}
Where the supremum is taken over all metrics $\sigma'$ conformally 
equivalent to
$\sigma$. (See Ahlfors \cite{ahlfors:invariants}.)
In particular, 
$\lambda_\sigma$ only depends on the conformal class of $\sigma$.

If $A$ is an annulus with conformal structure $[\sigma]$ then its {\em 
conformal modulus} $m(A,\sigma)$ (or $m(A)$ if $\sigma$ is implicit) is the 
extremal length of the class of arcs joining its boundaries. The
extremal length of the class of simple closed curves homotopic to
either boundary is $1/m(A)$.  A standard Euclidean metric on
$A$ is unique in its conformal class up to constant multiple, and 
realizes the supremum in the definition of the modulus. Thus $m(A) = H/L$ 
where $H$ is the distance between the boundaries in this metric, and $L$ 
is the length of either boundary.

\subsection{The Thick-Thin decomposition.}
\label{thick-thin}
Fixing a standard hyperbolic metric $\sigma$ on $S$, 
the {\em $\ep$-thin part} of $(S,\sigma)$
is the subset of $S$ where the injectivity radius is at most
$\ep$. The closure of its complement is called the $\ep$-thick part. 
If $\ep$ is less than the Margulis constant for $\Hyp^2$ (see e.g.
\cite{wpt:textbook}),
then the $\ep$-thin part  consists of compact annular neighborhoods of short
geodesics, and non-compact cuspidal components, or neighborhoods of 
punctures. Let us call a compact annulus $A$ of the thin part {\em 
peripheral} if its core is homotopic to $\boundary S$ (in which case one 
of its boundary components lies on $\boundary S$), and {\em internal} 
otherwise. In either case call a boundary component of $A$ {\em internal} 
if it lies in the interior of $S$. 

For our later convenience we shall describe the following slight 
variation on the  standard decomposition, which we shall call the {\em 
$(\ep_0,\ep_1)$  collar decomposition}, for two numbers
$0<\ep_1<\ep_0$ less than the Margulis constant.
Let $A_1,...,A_k$ denote the
annular neighborhoods of the geodesics of length
$\delta_i\le\ep_1$, such that the internal boundary components of $A_i$ have 
length $\ep_0$. 
It is not hard to see that the $A_i$ are
contained in separate components of the $\ep_0$-thin part, and that
$S-\union A_i$, minus neighborhoods of the cusps,  are in the
$\ep'_1$-thick part, for some $\ep'_1$ possibly slightly  smaller than $\ep_1$.

Let $\AAA$ denote the set of $A_i$, and let $\PP$ denote the set of 
closures $P$ of components of $S-\union A_i$. We call $\PP$ the set of 
{\em hyperbolic} components of the decomposition. 
Let $\QQ = \PP\union \AAA$ denote the entire decomposition.

It is well-known that 
$\ep_0$ may be chosen sufficiently small (independently of the surface)  
that, for any simple geodesic $\gamma$ in $S$, the components of 
$\gamma\intersect A_j$ are either segments that traverse $A_j$ from one 
boundary to the other, or are equal to the geodesic core of $A_j$.

One can see from an explicit computation (in the upper half-plane, for
example) that the modulus of each $A_i$ is 
\begin{equation}
\label{modulus for short geodesic}
	m(A_i) = {\pi\over\delta_i} - {2\over\ep_0},
\end{equation}
and we can guarantee by choosing $\ep_1$ sufficiently small
that $m(A_i) $ is larger than any given constant $m_0$. 

\subsection{The Teichm\"uller metric, and Kerckhoff's theorem.}
\label{kerckhoff theorem}
The {\em Teichm\"uller distance} on $\TT(S)$ is defined for two
conformal structures $[\sigma],[\tau]$ as
\begin{equation}
\label{define teich dist}
d_{\TT(S)}(\sigma,\tau)  = \half\log K(\sigma,\tau)
\end{equation}
where $K\ge 1$ is the least number such that there exists a
homeomorphism  homotopic 
to the identity on $S$ which is a $K$-quasiconformal map between the 
$\sigma$ and the $\tau$ conformal structures. We shall be using a 
different characterization of $K(\sigma,\tau)$, due to Kerckhoff, namely:

\begin{theorem}{Kerckhoff distance thm}{}
For a surface of finite topological type $S$, the following holds: 
\begin{equation}
K(\sigma,\tau) = \sup_{\alpha\in\CC(S,\boundary S) }
{\lambda_{S,\sigma}(\alpha)\over \lambda_{S,\tau}(\alpha)}
\label{Kerckhoff distance}
\end{equation}
\end{theorem}

\bfheading{Remarks.} Kerckhoff proved this in \cite{kerckhoff} for a 
closed surface $S$. However the same techniques apply to surfaces of 
finite analytic type, and the case with boundary can also be deduced 
easily using a doubling argument (see e.g. Abikoff
\cite[\S II.1.5]{abikoff} for an example of such arguments). 
An interesting feature of the
equivalence of the two definitions of $K$ is that the first one is fairly
easily seen to be symmetric in $\sigma$ and $\tau$ (since the inverse of 
a $K$-quasiconformal map is also $K$-quasiconformal), but the symmetry of 
the second is less evident. We shall use this to advantage in section 
\ref{product structure}.

\medskip

The Teichm\"uller space of the torus is isometric to the hyperbolic
plane $\Hyp^2$, and Kerckhoff's theorem in this case can be written as
a simple formula, after proper interpretation. We will use this
formula not for the torus but in conjunction with twisting numbers and
the annulus contributions to extremal length, $\lambda_{A,\sigma}$
(definition (\ref{A length def})) in section \ref{product
structure}. For the sake of exposition we shall give an elementary
proof:

\begin{lemma}{Kerckhoff-Royden for torus}{}
Let $\Hyp$ denote the upper half-plane $\{z: Im(z)>0\}$, and define for
$z_j = x_j+iy_j$
$$
K(z_1,z_2)  = \sup_{t\in\R} { y_2 + (t+x_2)^2/y_2 \over
				     y_1 + (t+x_1)^2/y_1}.
$$
Then
\begin{equation}
\label{hyp distance is log K}
d_\Hyp(z_1,z_2) = \half\log K(z_1,z_2)
\end{equation}
where $d_\Hyp(z_1,z_2)$ is the hyperbolic distance function.
\end{lemma}
(Note: actually the common definition of $d_\Hyp$ does not
include the factor of a half; however we will retain it in order to be
consistent with (\ref{define teich dist}).)

\begin{proof}{}
Consider first the case that $z_1$ and $z_2$ are purely imaginary.  It
is an easy exercise that $K(iy_1,iy_2) = \max\{y_2/y_1,y_1/y_2\}$.
Since $\half|\log(y_2/y_1)|$ is exactly the hyperbolic distance
function, this completes the proof in this case. 

Since for general $z_1,z_2\in \Hyp$ there is always some hyperbolic isometry
$A\in \SL 2(\R)$ such that $A(z_1),A(z_2)$ are purely imaginary, it
would suffice to show that $K$ is invariant under the action of $\SL
2(\R)$. 

One can see this conveniently by reducing it to a linear algebra
problem. Let $\widehat
\Hyp$ denote $\{(\alpha,\beta)\in\C^2: Im(\beta/\alpha)>0\}$. For
$(\alpha,\beta)\in\widehat\Hyp$ and $u,v\in\R$, define
$$
\lambda_{\alpha,\beta}(u,v) = { |u\alpha + v\beta|^2 \over
					Im(\beta\overline\alpha)}.
$$
Then define $\widehat K$ on $\widehat \Hyp\times\widehat\Hyp$ by
$$
\widehat K(\alpha_1,\beta_1,\alpha_2,\beta_2)  = 
\sup_{u,v\in\R} { \lambda_{\alpha_2,\beta_2}(u,v) \over
                  \lambda_{\alpha_1,\beta_1}(u,v)}.
$$
One can check that the left action of $\SL 2(\R)$ on $\C^2$ (viewed as
column vectors)
preserves $Im(\beta\bar\alpha) = |\alpha|^2Im(\beta/\alpha)$, and
therefore in particular preserves 
$\widehat\Hyp$.
Furthermore note that the denominator of $\lambda_{\alpha,\beta}(u,v)$
is the squared norm of an inner product $(u,v)\cdot(\alpha,\beta)$, so
it is easy to see 
that $\lambda_{A(\alpha,\beta)}(u,v) = \lambda_{\alpha,\beta}((u,v)A)$.
It follows that $A\in \SL 2(\R)$ 
leaves $\widehat K$ invariant. 
since the supremal ratio over all $(u,v)$ is the same as the supremum
over all $(u,v)A$.

It is an easy computation to show that $\lambda_{1,z}(t,1) = y +
(t+x)^2/y$ if $z=x+iy$, and that $\lambda_{1,\beta/\alpha}(u,v) =
\lambda_{\alpha,\beta}(u,v)$. Moreover, $\lambda_{\alpha,\beta}(cu,cv)
= c^2\lambda_{\alpha,\beta}(u,v)$. 
It follows that 
$$
\widehat K(\alpha_1,\beta_1,\alpha_2,\beta_2)  =
K(\beta_1/\alpha_1,\beta_2/\alpha_2). 
$$
The map $(\alpha,\beta)\mapsto \beta/\alpha$ is a fibration of
$\widehat\Hyp$ over $\Hyp$, which takes the linear action 
of $\SL 2(\R)$ upstairs to its action by hyperbolic isometries downstairs. 
We conclude that $K$ is also invariant under $\SL 2(\R)$.  
\end{proof}

\bfheading{Remark.}
One can give
a geometric interpretation to this proof by associating to 
$(\alpha,\beta)\in\widehat\Hyp$ a torus, obtained by identifying
opposite sides on a parallelogram spanned by $\alpha$ and $\beta$ in
$\C$. The quantity $Im(\beta\overline\alpha)$ is the area of this
parallelogram. For $u,v\in\Z$, the quantity
$\lambda_{\alpha,\beta}(u,v)$ is the extremal length of the homotopy class
in the torus represented by $(u,v)$ in the generators $\alpha$ and
$\beta$ (for general $u,v\in\R$, one can give an interpretation in
terms of measured foliations on the torus). Now $\widehat K$ can be
seen as a supremum of  extremal length ratios, just as in Kerckhoff's
theorem. The fibration from $\widehat \Hyp$ to $\Hyp$ is just
identification under conformal equivalence.

\section{Twisting numbers and Fenchel-Nielsen coordinates}
\label{twisting section}
In this section we will describe a measure, called
$t_{\gamma,\sigma}(\alpha)$, for
the amount that a simple curve $\alpha$ in a hyperbolic surface $(S,\sigma)$
``twists'' around another curve $\gamma$. In particular if $\gamma$ is
short then this number will approximate the amount of twisting that
$\alpha$ undergoes in the thin annulus corresponding to $\gamma$, and
this will determine the contribution of this annulus to the extremal
length of $\alpha$. 
At the end of the section we will relate this construction to the
familiar Fenchel-Nielsen twist coordinates.

Let $\sigma$ be a standard hyperbolic metric on $S$. Let $\gamma$
denote a (homotopically non-trivial) oriented simple closed curve, and
let $\alpha$ denote any element of $\CC'(S,\boundary S)$, which
intersects $\gamma$.  Let $x$ be a
point of the intersection $\gamma^\sigma\intersect \alpha^\sigma$. 
In the universal cover $\til S$ of $S$, which is a convex subset of $\Hyp^2$,
let $L_\gamma$ be a geodesic in the lift of $\gamma^\sigma$, and
$L_\alpha$ be a geodesic in the lift of $\alpha^\sigma$ (extended to
an infinite geodesic if it meets $\boundary \til S$), such that
$L_\gamma$ and $L_\alpha$ meet in a point that projects to $x$. 
Identify $L_\gamma$ isometrically with $\R$, consistent with the orientation
of $\gamma$, and let  
let $a_R,a_L$ be 
the endpoints of $L_\alpha$ to the right and left, respectively, of $L_\gamma$.
Let $p_\gamma:\Hyp^2\union S^1\to L_\gamma$ be the orthogonal 
projection, and define the (signed) twisting number
\begin{equation}
\label{twisting number def}
t'_{\gamma,\sigma}(\alpha,x) = {p_\gamma(a_R) - p_\gamma(a_L)\over 
\ell_\sigma(\gamma)}.
\end{equation}
(see figure \ref{twisting def fig}). 
It is clear that the definition does not depend on the particular
choice of lifts. 

\realfig{twisting def fig}{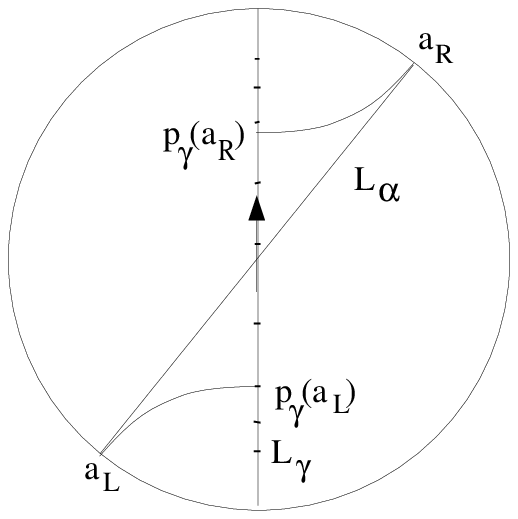}{}

\begin{lemma}{twist for conjugates}{}
Let $\gamma$ be a homotopically non-trivial oriented simple closed
curve, $\alpha\in\CC(S,\boundary S)$, and $\sigma$ a standard hyperbolic metric on
$S$. If $x$ and $x'$ are two points of $\gamma^\sigma \intersect
\alpha^\sigma$, then 
$$|t'_{\gamma,\sigma}(\alpha,x)-t'_{\gamma,\sigma}(\alpha,x')| \le 1.$$
\end{lemma}
\begin{proof}{}
Let $g$ be the indivisible element of $\pi_1(S)$ whose axis is $L_\gamma$. 
Choosing appropriate lifts $L_\alpha$ and $L'_\alpha$ corresponding to
$x$ and $x'$, we may assume that
$L'_{\alpha}$ is trapped between 
$L_\alpha$ and $g(L_\alpha)$. The proof is completed by figure 
\ref{twist for conjugates fig}.
\end{proof}

\realfig{twist for conjugates fig}{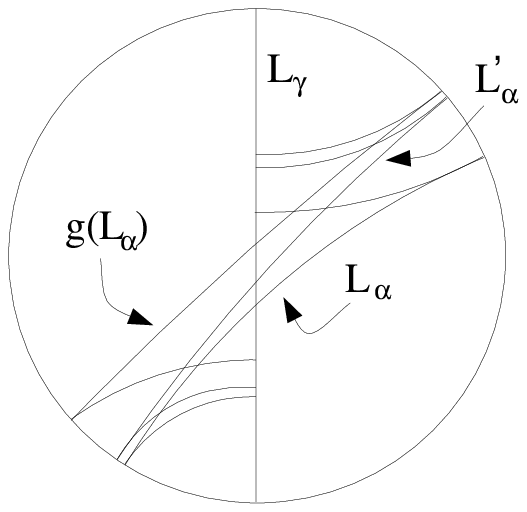}{}

Thus it is useful to define the quantity
\begin{equation}
\label{inf twist def}
t_{\gamma,\sigma}(\alpha) = \min_{x} t'_{\gamma,\sigma}(\alpha,x)
\end{equation}
where the infimum is taken over all $x$ in $\gamma^\sigma\intersect
\alpha^\sigma$. 

We can also measure twisting inside an annulus. If $A$ is an annulus with 
conformal structure $[\sigma]$ then let $\sigma$ be a standard 
Euclidean metric on $A$ in this conformal class. The univeral cover of 
$A$ is isometric to a strip $[0,L]\times\R$ in the plane, with 
fundamental domain of height $H$. If $\alpha$ is a curve
connecting the  boundaries of $A$, it lifts to a curve with endpoints 
$(0,y_0)$ and $(L,y_1)$, and we may define 
$$t_{A,\sigma}(\alpha) = 
(y_1-y_0)/H$$ 
(the sign of $t$ is again determined after a choice of 
orientation of $A$). 

Now consider again a standard hyperbolic metric $\sigma$ on $S$, and an 
annulus $A$ of the $(\ep_0,\ep_1)$ collar decomposition of $(S,\sigma)$. 
If $\alpha\in\CC(S)$ then its geodesic 
representative $\alpha^\sigma$ intersects $A$
in several components, whose twists differ by no more than one. 
Thus we may take $t_{A,\sigma}(\alpha)$ to be the minimum of these. 
Note in particular that if $A$ is peripheral then 
$t_{A,\sigma}(\alpha) = 0$ since $\alpha^\sigma$ is perpendicular 
to $\boundary A$. 

\begin{lemma}{internal and external twists}{}
There exist $\ep_0,\ep_1$ less than the Margulis constant such that, 
for any standard hyperbolic metric $\sigma$ on a surface $S$,
if $A$ is an $(\ep_0,\ep_1)$ collar with (non-peripheral) core $\gamma$ 
and $\alpha$ is any element of $\CC(S)$ that crosses $\gamma$, 
$$
    |t_{A,\sigma}(\alpha) - t_{\gamma,\sigma}(\alpha)| \le 1.
$$
\end{lemma}

\begin{proof}{} 
In case $A$ is peripheral both twists are 0, so let us assume
$A$ is internal.
Consider figure \ref{internal-external fig}, where in the 
upper half-plane model of $\Hyp^2$ we place a lift $L_\gamma$ of 
$\gamma^\sigma$, and the corresponding lift $\til A$ of $A$. 
Pick $x\in\gamma^\sigma\intersect \alpha^\sigma$ and choose a
component $L_\alpha$ of the lift of $\alpha^\sigma$ which intersects
$L_\gamma$ in a point projecting to $x$. Let $\alpha_x$ denote the
component of $\alpha^\sigma\intersect A$ which contains $x$.
We may arrange 
the figure so that $L_\alpha$ exits $\til A$ as shown, terminating at infinity 
in the interval $J$. 

\realfig{internal-external fig}{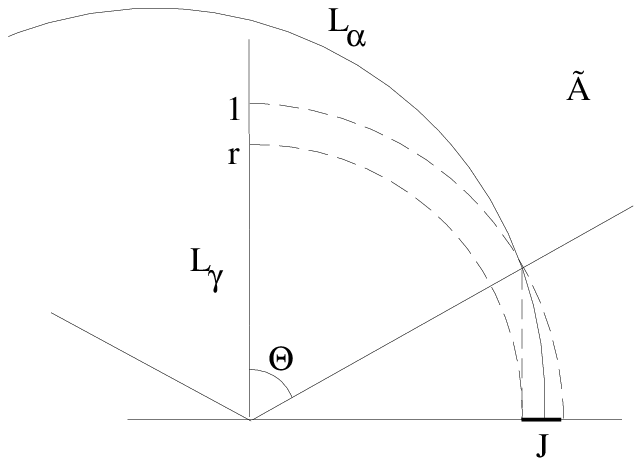}{}

The angle $\theta$ between $L_\gamma$ and $\boundary \til A$ satisfies $\cos 
\theta \le \ep_1/\ep_0$, and the position $r$ of the projection of the 
endpoint of $L_\alpha$ satisfies $1\ge r > \sin \theta$. It follows 
immediately that 
$$
0\le t'_{\gamma,\sigma}(\alpha,x) - t_{A,\sigma}(\alpha_x) < {1\over \ep_1} \left| \log 
(1-\ep_1^2/\ep_0^2)\right|.
$$
Thus given $\ep_0$  we may choose $\ep_1$ sufficiently small that this is 
at most 1. Since both
$t'_{\gamma,\sigma}(\alpha,x)$ and $ t_{A,\sigma}(\alpha_x)$
are minimized for the same $x$ (both depend monotonically on the angle
at which $\gamma^\sigma$ and $\alpha_x$ cross), we have the desired inequality.
\end{proof}

The following lemma observes that for a 
thin annulus $A$ any twisting number can be achieved, up to an additive 
constant, by some curve that spends ``most'' of its time in $A$. This fact 
will be used in the proof of theorem \ref{Product structure}. 

\begin{lemma}{closed curve from annulus segment}{}
Let $(Q,\tau)$ be a surface of finite topological type with standard 
hyperbolic metric $\tau$.
Let $\ep_0 >\ep_1>0$ be as in lemma \ref{internal and external twists}, 
and fix $\ell_1>0$. 
Suppose that $\ell_\tau(\gamma) \le \ell_1$ for all boundary 
components $\gamma$ of $Q$, 
and let $A\subset Q$ be an internal annulus of the $(\ep_0,\ep_1)$ collar
decomposition.
Given $t\in\R$ there exists $\alpha\in\CC(Q)$ such that 
$$|t_{A,\tau}(\alpha)-t|\le c_1,$$
and $\alpha$ consists of one or two arcs
traversing $A$, together with one or two segments of hyperbolic length
at most $c_2$.
The constants $c_1,c_2$ depend only on $\ep_0,\ep_1,\ell_1$ and
$\chi(Q)$.
\end{lemma}

\begin{proof}{}
We claim that, for some $r_1$ depending on $\ep_0,\ep_1,
\ell_1$ and $\chi(Q)$,  there is, for each boundary component of $A$, 
a geodesic arc outside $A$ of length at most $r_1$, which either joins 
that component to the other, or to itself. 
To see this, let $A_r$ denote the $r$-neighborhood of $A$ in 
$Q$; $A_r$ is the union of $A$ with two pieces, each adjacent to a 
boundary component. 
If at least one of them is a collar, i.e.
an embedded annulus of radius $r$ in the interior of $Q$,   we have 
$\Area(A_r) >  
\ep_1e^r$, and since $\Area(Q) \le 2\pi|\chi(Q)|$ we have an upper bound 
$r< r_0=\log 2\pi|\chi(Q)|/\ep_1$. 
Thus, if $r\ge r_0$ then both added pieces fail to be collars. 
If for some $r\le r_0$ a collar meets itself or the other collar we have
the desired geodesic arc of length $r_1 = 2r_0$. 	The other possibility 
is for a collar to hit the boundary at some $r\le r_0$. 
In this case there is an arc of length at most $r_1 = 2r_0 + \ell_1$ 
which runs once around the boundary and back to $A$. 

Now if there is an arc joining both components, 
we can take $\alpha_0$ to be the union of this short arc 
with an arc that runs through $A$. If there are two 
arcs, joining each boundary to itself, we can join these arcs to a union 
of two parallel arcs through $A$. Applying $n$ Dehn twists around the
core of $A$ to $\alpha_0$ will clearly increase (or decrease) its twisting
number by about $n$ (see also lemma \ref{compare twists}), so that by
an approriate twist we may obtain 
$\alpha$ with the desired twisting number.
\end{proof}

In the proof of corollary \ref{generalized estimate},
we will need to measure twisting using a hyperbolic 
metric that is not quite standard. In particular, suppose that $\sigma$ 
is a standard hyperbolic metric on $S$ and that $P$ is a 
subsurface consisting of a union of components of the $(\ep_0,\ep_1)$
collar decomposition, such that any annular components in $P$ are
non-peripheral in $P$. 
The metric $\sigma$ restricted to $P$ is not 
standard, but it is  
not too far from it. Let $\sigma'$ be the standard hyperbolic metric on 
$P$ which is conformally equivalent to $\sigma$. Now if $\gamma\in P$ is 
a non-peripheral geodesic of length $\ell_\sigma(\gamma)<\ep_1$, we can 
measure twisting around $\gamma$ with respect to both metrics. The 
following lemma compares the two: 

\begin{lemma}{twist in subsurface}
Let $P, S,  \sigma$ and $\sigma'$ be as above. There exist $\ep_0,\ep_1$ 
sufficiently small and $c>0$ such that, if 
$\ell_\sigma(\gamma^\sigma)<\ep_1$ for a  
non-peripheral simple closed curve $\gamma\subset P$, then
$$
|t_{\gamma,\sigma}(\beta) - t_{\gamma,\sigma'}(\RR_P(\beta)) 
| \le c
$$
for any $\beta\in\CC(S,\boundary S)$.
\end{lemma}

\begin{proof}{}
We shall give the proof in the case that $S$ has no boundary. The
general case can be reduced to this by means of a doubling argument.
Let $\beta^\sigma$ be the $\sigma$-geodesic representative of $\beta$. 
Let $P_0$ be the complement in $P$ of $\ep_\MM$-thin annulus
neighborhoods 
of $\boundary P$ (measured in $\hat\sigma$), where $\ep_\MM>\ep_0$ is 
the Margulis constant.  
We claim that, given $\nu>0$ there is a choice of $\ep_0$ so that in
$P_0$ the ratio $\sigma/\hat\sigma$ is $\nu$-close to 1, and its derivatives 
(in normal $\hat\sigma$ coordinates) are bounded by $\nu$. This will
imply that  
$\beta^\sigma\intersect P_0$ is a curve system with very low ($O(\nu)$) 
curvature in the metric $\hat\sigma$ . Thus the $\hat\sigma$-geodesic
$\hat\beta_0$ homotopic to  
$\beta^\sigma\intersect P_0$ rel endpoints is $O(\nu)$-close to it. 
Now, we are interested in $\RR_P(\beta)^{\hat\sigma}$, the 
$\hat\sigma$-geodesic  representative of $\beta^\sigma\intersect P$, rel 
$\boundary P$. But by lemma \ref{bounded slide} applied in $(P,\hat\sigma)$,
the intersection $\RR_P(\beta)^{\hat\sigma}\intersect P_0$  is obtained by a 
bounded homotopy from $\hat\beta_0$. It follows that the 
twisting number around $\gamma$, which is certainly in $P_0$, has not 
changed by more than a bounded amount.

It remains to prove our claim about $\sigma/\hat\sigma$, but this is 
an elementary argument in complex analysis.
If we consider the universal coverings of $P$ in (subsets 
of) the unit disk which are obtained from the two metrics $\sigma$ and 
$\hat\sigma$, we see that the identity on $P$ lifts to a conformal map $\Phi$
between these covers. For $x\in P_0$, the hyperbolic distance to 
$\boundary P$, in either metric, is roughly $\log \ep_\MM/\ep_0$, which 
we can make as large as we please. If $x$ lifts to $0$ in both covers we 
see (using the formulas for hyperbolic distance) that $\Phi$ must map a 
disk of radius roughly $(\ep_\MM-\ep_0)/(\ep_\MM+\ep_0)$ univalently over 
a disk of similar radius.  The Schwarz lemma applied to $\Phi$ and
$\Phi^{-1}$ 
implies that $|\Phi'(z)|$ is close to 1 in a neighborhood of 0, and an
application of the Cauchy integral formula serves to bound
$\Phi''(0)$. Our claim about $\sigma/\hat\sigma$ follows. 
\end{proof}

%
%
%

Recall now the construction of 
Fenchel-Nielsen coordinates on a Teichm\"uller space $\TT(S)$ (see e.g.
\cite{abikoff}) 
Let $\gamma_1,\ldots,\gamma_n$ denote a system of oriented simple closed curves
which cuts $S$ into a union of pairs of pants. Suppose, to start, that
$S$ has no punctures and $p$ boundary components, so that
$n=-3\chi(S)/2$ and each pair of pants is compact. 
For any standard hyperbolic structure $\sigma$ on $S$ we
immediately obtain $n$ positive numbers
$\{\ell_\sigma(\gamma_i^\sigma)\}$, the {\em length parameters.}
Another $n-p$ real numbers describe ``twisting'' parameters for the
gluings of the pairs of pants. Let us make this concrete as follows. 

A pair of pants $P$ with boundary curves
$\gamma_1,\gamma_2,\gamma_3$ 
contains three unique homotopy classes of simple properly embedded arcs
$\alpha_{12},\alpha_{23},\alpha_{13}$ (the ``seams''), such that
$\alpha_{ij}$ joins $\gamma_i$ to
$\gamma_j$. 
Fix a set of representatives of the seams which match on opposite
sides of each non-peripheral $\gamma_i$. This determines a system of
curves $\mu$ in $\CC(S,\boundary S)$. We call the (oriented) $\gamma_i$
together with $\mu$ a {\em marking}, and note that it is purely
topological information.

Now for any standard hyperbolic structure $\sigma$ on $S$, the seams
have geodesic representatives $\alpha_{ij}^\sigma$ which are orthonal
to $\gamma_i^\sigma$ and $\gamma_j^\sigma$. 
The seams cut each pair of pants into two congruent right-angled
hexagons, and in particular they bisect the boundary components. Now
for each $\gamma_j$ and each pair of corresponding seam endpoints
there is a unique geodesic path along
$\gamma_j$ that must be spliced between the endpoints,
so that the resulting  curve is in the homotopy class $\mu$. 
The length $m_j(\sigma)$ of this path (which can be given a sign,
since $\gamma_j$ is oriented) is the same for both pairs of seam
endpoints on $\gamma_j$. 
We define our {\em twist parameter}  $s_j(\sigma)$ to be $m_j(\sigma)
/ \ell_\sigma(\gamma_j)$. 

This construction gives a homeomorphism $F: \TT(S) \to
\R_+^n\times\R^{n-p}$, taking $\sigma$ to 
$(\ell_\sigma(\gamma_1),\ldots,\ell_\sigma(\gamma_n),
s_1(\sigma),\ldots,s_{n-p}(\sigma))$.
We note also that a positive Dehn twist on
$\gamma_j$ has the effect of incrementing $s_j(\sigma)$ by one and
leaving the other coordinates invariant.

In the case with punctures the construction is similar, except that
one or two  ends of a pair of pants may be a puncture, and we allow
the seams to be non-compact arcs terminating in boundaries or punctures, 
so that $\mu$ is in $\CC'(S,\boundary S)$.

\medskip

To conclude this section, we give a comparison between our
twisting numbers and the Fenchel-Nielsen twist coordinates. 

\begin{lemma}{compare twists}
Let $\gamma_1,\ldots,\gamma_n$ and $\mu$ be a marking for $S$, and
$s_j:\TT(S)\to \R$ be the associated twist coordinates. 
Then for any standard hyperbolic metric $\sigma$ on $S$,
$$
|s_j(\sigma) - t_{\gamma_j,\sigma}(\mu)| \le 1.
$$
Furthermore, given any  $\alpha\in\CC(S,\boundary S)$ and two standard
hyperbolic metrics $\sigma$ and $\sigma'$,  
$$
\big| (t_{\gamma_j,\sigma}(\alpha) - t_{\gamma_j,\sigma'}(\alpha))
 - (s_j(\sigma) - s_j(\sigma')) \big| \le 4
$$
\end{lemma}

\begin{proof}{}
To obtain the first inequality, lift $\gamma_j^\sigma$ to a geodesic $L_j$ in
$\Hyp^2$, and consider a lift $\til\mu$ of a component of $\mu^\sigma$ that
intersects $L_j$. This lift is homotopic to a curve that follows two
(lifts of) seam curves $\alpha_{jk}^\sigma,\alpha_{jm}^\sigma$ whose
endpoints on $L_j$ are separated by 
$s_j(\sigma)$. The other
endpoints of these seams meet lifts $L_k$ and $L_m$ of
$\gamma^\sigma_k,\gamma^\sigma_m$. 
Sinc $L_k$ (resp. $L_m$) is disjoint from its translates, its
projection to $L_j$ has length at most $\ell_\sigma(\gamma_j)$. The
endpoints of $\til \mu$ at infinity must be separated from $L_j$ by $L_k$
and $L_m$, respectively. It follows that the length of the projection
of $\mu$ to $L_j$, divided by $\ell_\sigma(\gamma_j)$,  differs from
$s_j(\sigma)$ by at most 1.


For the rest of the lemma we note the fact that, for any two curve
systems $\alpha,\beta \in\CC'(S,\boundary S)$, the {\em difference}
$t_{\gamma_j,\sigma}(\alpha) - t_{\gamma_j,\sigma}(\beta)$ is
independent of $\sigma$ up to a bounded error (of 1). This is because
we can detect this difference topologically:
lifting to the universal cover, let $L_\alpha$ be a lift of a
component of $\alpha$ which crosses a lift $L_j$ of $\gamma_j$. The
translates of $L_\alpha$ along $L_j$ divide the plane into a sequence
of consecutive strips, and a component of a lift of $\beta$ that
crosses $L_j$ meets
some number $\delta$ of these strips. (The number obtained from a
different lift of $\beta$ may differ by at most 1). It is easy to see that
$\delta$, appropriately signed, gives this twist
difference, up to an error of 1.  

The second inequality now follows from this fact applied to $\alpha$
and $\mu$, together with the first inequality.
\end{proof}
\section{Basic extremal length estimates}
\label{basic extremal}
In this section we shall discuss and derive some fairly well-known 
estimates for extremal lengths, which will be the building blocks of the 
proof in section \ref{main length estimate}. 

Here and in the rest of
the paper, we shall use ``$a\simmult b$'' to mean that the ratio $a/b$
is bounded above and below by positive constants, which usually depend
on topological data, and on previously chosen constants.
Similarly, ``$a\simle b$'' and ``$b\simge a$'' mean that $a/b$ is
bounded above. Where necessary we
will indicate explicitly the dependence of the constants.

\subsection{Spacing and upper bounds.}
Lower bounds on extremal lengths are relatively easy to obtain, by
exhibiting some metric and computing the length in that metric. The
following lemma is helpful for obtaining upper bounds. 

For a representative curve $\hat\beta$ of a class $\beta$ in
$C(S,\boundary S)$ and a metric $\sigma$ on $S$, define the {\em
spacing} $\nu_{S,\sigma}(\hat\beta)$ as the supremum of numbers $\nu$ such that
the $\nu$-neighborhood $\NN_\nu(\hat\beta)$ is a standard product
$\hat\beta\times[0,1]$. We have:

\begin{lemma}{spacing estimate}{}
Given a metric $\sigma$ on $S$ and a representative $\hat\beta$ of $\beta\in
C(S,\boundary S)$, the extremal length of $\beta$ is bounded by
$$
	\lambda_{S,\sigma}(\beta) \le {\Area(S,\sigma)\over
					4\nu^2_{S,\sigma}(\hat\beta)}.
$$
\end{lemma}

\begin{proof}{}
Let $\beta_1,\ldots,\beta_n$ denote the components of $\hat\beta$. The
neighborhoods $R_i=\NN_\nu(\beta_i)$ for any 
$\nu<\nu_{S,\sigma}(\hat\beta)$ are, by assumption, either annuli
or rectangles. In either case, the extremal length of the family of
transverse arcs $\alpha_i$ connecting the boundaries
$\beta_i\times\{0,1\}$ is, applying the definition,
\begin{equation}
	\lambda_{R_i,\sigma}(\alpha_i) \ge {4\nu^2\over\Area(R_i,\sigma)}.
\label{alpha bound}
\end{equation}
Let $[\beta_i]$ denote the family of curves in $R_i$ freely homotopic to 
$\beta_i$ when $\beta_i$ is a closed curve, and homotopic with endpoints 
on the edges $\boundary\beta_i\times[0,1]$ when $\beta_i$ is a segment. 
Then it is a standard fact that 
$$\lambda_{R_i,\sigma}([\beta_i]) = 1/\lambda_{R_i,\sigma}(\alpha_i).$$
This together with (\ref{alpha bound}) 
bounds $\lambda_{R_i,\sigma}([\beta_i])$ from above. 
Now let $\beta_R$ denote the class of all curve systems that are a union 
of representatives from each $[\beta_i]$. 
Since the $R_i$
are all disjoint it is easy to see that 
\begin{equation*}
\lambda_{\union R_i,\sigma}(\beta_R) = \sum_{i=1}^n 
\lambda_{R_i,\sigma}([\beta_i]) 
\end{equation*}
by an application of the Schwarz inequality (actually only the
``$\le$'' direction is needed, but the other direction is even easier).

The proof is completed by a standard monotonicity argument:
Since the class $\beta$ in $\CC(S,\boundary S)$ is bigger than the class 
$\beta_R$ (it has representatives not limited to the $R_i$), its 
length in any metric is smaller. 
Furthermore the area of $S$ in any metric is at least that of $\union R_i$.
Thus $\lambda_{S,\sigma}(\beta)\le \lambda_{\union R_i,\sigma}(\beta_R)$, 
and the lemma follows.
\end{proof}

\subsection{Extremal length in hyperbolic subsurfaces.}

If a hyperbolic subsurface $M$ of $S$ is ``insulated" from $S$ by sufficiently 
thick annuli on its boundary, then extremal lengths in $\CC(S)$ and 
$\CC_0(M)$ (curves  that are non-peripheral in $M$) 
are  approximately the same:
\begin{lemma}{insulated surface}{\normalshape \cite[lemma 8.4]{minsky:2d}}
Let $(S,\sigma)$ be any Riemann surface. There exist $m_0$ and $C$,
depending only on the topological type of $S$, such that
if  $M\subset S$ is a hyperbolic subsurface for which each component
$\gamma$ of  
$\boundary M$ bounds an annulus $A_\gamma\subset M$ with modulus
$m(A_\gamma)\ge m_0$, then for any
$\alpha\in\CC_0(M)$,
$$
1 \le {    \lambda_{M,\sigma}(\alpha)\over
        \lambda_{S,\sigma}(\alpha)} \le C.
$$
\end{lemma}

More generally, a curve in $S$ may cut through $M$, and if 
$M$ is a nicely chosen subsurface then we can say something about the 
extremal length of the intersection. For a standard hyperbolic metric 
$\sigma$, 
recall the $(\ep_0,\ep_1)$ collar decomposition $\QQ = \PP\union\AAA$.
For a hyperbolic component $P\in\PP$, define
\begin{equation}
\label{P length def}
\lambda_{P,\sigma} (\beta) = \lambda_{P,\sigma}(\RR_{P}(\beta)).
\end{equation}
Note that, provided $\ep_0$ is chosen sufficiently small, 
$\RR_{P}(\beta)$ is represented by $\beta^\sigma\intersect P$. This 
hyperbolic geodesic also gives a good estimate for extremal length on 
$P$, via the following:

\begin{lemma}{extremal is hyperbolic for thick}{}
Suppose that $\sigma$ is a standard hyperbolic metric for $S$, and
suppose that 
$\ell_\sigma(\gamma)\le \ell_0$ for any component $\gamma$ of $\boundary 
S$. 
For any $\beta\in\CC(S,\boundary S)$ and a hyperbolic component $P$ of the 
$(\ep_0,\ep_1)$ collar decomposition, 
$$
\lambda_{P,\sigma} (\beta) \simmult 
\ell^2_\sigma(\beta^\sigma\intersect P)
$$
Where the multiplicative constants for $\simmult$ depend only on the 
topological type of $S$ and on $\ep_0,\ep_1$ and $\ell_0$. 
\end{lemma}

\begin{proof}{}
Because of the absence of thin parts in $P$ we can rearrange 
$\beta^\sigma\intersect P$ to have sufficiently large spacing in the 
hyperbolic metric $\sigma$. Indeed,

\begin{lemma}{adjust on P}
There is a curve system $\beta'$ in $P$ homotopic to 
$\beta^\sigma\intersect P$ by a homotopy that moves no point of 
$\beta^\sigma\intersect \boundary P$ more than halfway around the 
boundary, and such that 
$$
\nu_{P,\sigma}(\beta') \ge {c\over \ell_\sigma(\beta^\sigma\intersect 
P)}
$$
where the constant $c$ depends only on the topological type of $S$,
and on $\ep_0,\ep_1$ and $\ell_0$
\end{lemma}

\begin{proof}{}
Subdivide $P$ into pairs of pants $Q_1,\ldots,Q_p$, using the shortest 
collection of pairwise disjoint simple geodesics in $P$. The number $p$ 
of pieces is given by $p=|\chi(P)|$, and we note that some pairs of 
pants may have one or two punctures rather than boundaries.
The lengths of all boundary components of $Q_i$ are bounded below by
$\ep_1$, and above by
a constant $\ell_1$ depending on $\ell_0$ and the topological type of
$S$.  (This is a standard fact, whose proof is similar to the proof
that the injectivity radius of $(S,\sigma)$ is bounded above in terms
of the hyperbolic area of $S$).

Begin by performing a homotopy on $\beta^\sigma\intersect Q_j$, for
each $j$, which moves 
the points of $\beta^\sigma\intersect 
\boundary Q_j$ around
$\boundary Q_j$ until they are evenly spaced, and keeps the rest of the 
curve geodesic. Call the new curve $\beta_1$. Clearly this can 
be done so that no endpoint 
moves more than halfway around its boundary component. Furthermore, provided 
$\ep_0$ has been chosen sufficiently small, $\beta_1$ will still be 
contained in $P$. 

The bound $\ell_1$ on boundary lengths means that 
the length of no arc has changed by more than an additive constant by 
this homotopy. It also means that there is a lower bound $\ell_2$ on the 
length of any path joining two boundaries of any $Q_j$. Furthermore no 
closed component of $\beta^\sigma\intersect P$ has length less than
$\ep_1$. We conclude that 
$$
\ell_\sigma(\beta_1) \le c_1\ell_\sigma(\beta^\sigma\intersect P)
$$
where $c_1$ depends on the previous constants.

Now we deform $\beta_1$ in each $Q_k$ separately to obtain a curve 
$\beta_2$ with maximal 
spacing. Consider first the case where $\beta_1\intersect Q_k$ has no 
closed curve components, and where $Q_k$ has three boundaries and no 
cusps.

One way to make the construction concrete is
to form the following Voronoi diagram: For $i=1,2,3$ let $C_i$ be
the collar defined by 
$$C_i=\{x\in Q_k: d(x,\gamma_i) \le d(x,\gamma_j) \ \hbox{for}\ j\ne i\}.$$
Then each $C_i$ is an annulus whose whose width at any point is
bounded between two constants $c_2,c_3$, and whose
boundary in the interior of $Q_k$
consists of two geodesic arcs
$\alpha_{ij} = C_i\intersect C_j$ (for $j\ne i$) of length
$s_{ij}\in [c_4,c_5]$
(The estimates here are uniform since the lengths
of $\gamma_i$ are bounded above and below, so we have a compact set of
possible pairs of pants. Alternatively one can obtain explicit
estimates using hyperbolic trigonometry.)

In fact, one can use (for example) the foliation of each $C_i$ by
geodesics perpendicular to $\gamma_i$ to obtain a 
Euclidean metric $\hat\sigma$ on $C_i$ which is bilipschitz equivalent
to $\sigma$, with a uniform constant.

Now for each arc $\tau$ of $\beta_1\intersect Q_k$, associate a point
$x(\tau)$ on $\alpha_{ij}$ if the endpoints
of $\tau$ lie on $\gamma_i$ and $\gamma_j$. Arrange the points
$\{x(\tau)\}$ so that they are uniformly spaced on each $\alpha_{ij}$.
There is a unique arc $\tau'$, homotopic to $\tau$,  which spirals
in $C_i$ from $\gamma_i$ 
to $x(\tau)$, and then in $C_j$ from $x(\tau)$ to $\gamma_j$, and 
is a $\hat\sigma$-geodesic 
in each $C_i$. The homotopy from $\tau$ to $\tau'$ 
will move no point more than a bounded distance (this is easiest to
see in the universal cover -- see figure \ref{pants diagram}).
It is also easy to see that the resulting arcs will be disjoint

\realfig{pants diagram}{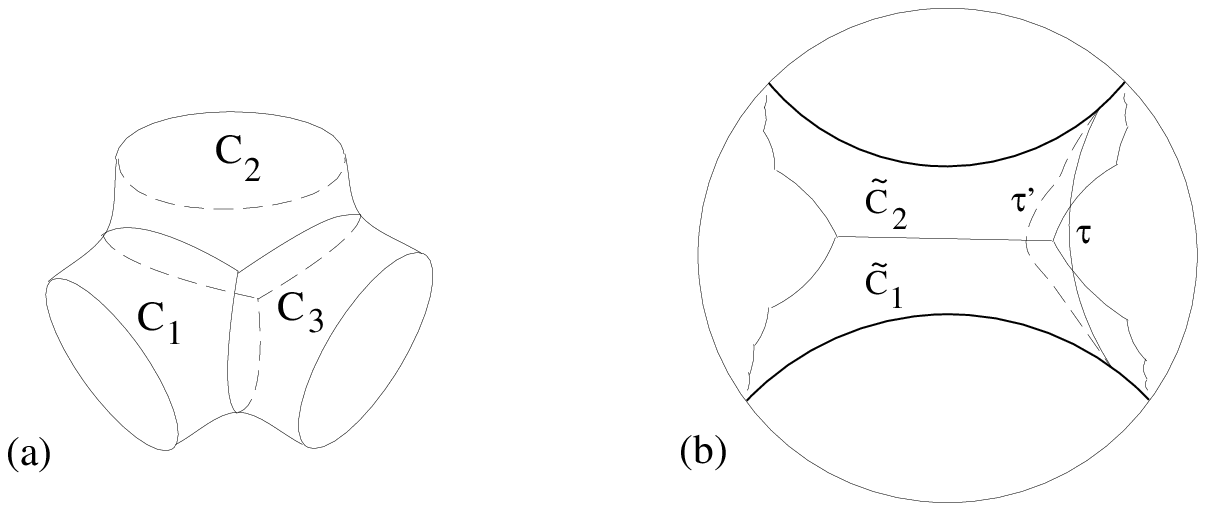}{(a) The decomposition of a 
pair of pants into  
three annuli. (b) In the lifts of two adjacent annuli to $\Hyp^2$, an
example of an arc $\tau$ and the homotopic arc $\tau'$.}

We can treat the case with cusps similarly, except that we first excise a 
horospherical neighborhood of each puncture with boundary length $\ep_0$, 
and work in the complement. 

If $\beta_1$ has a closed component in $Q_k$, it must be homotopic to a 
boundary component $\gamma_j$, and no part of $\beta_1$ can cross that 
component. Thus we can take a collar neighborhood of $\gamma_j$ of 
definite width and space the corresponding components evenly (and 
maximally) there. 

Now for maximally spaced curves in a Euclidean annulus it is clear that the 
spacing is approximately the area divided by the total length.
Where two collars meet along a common boundary segment, the local picture 
is, up to uniform bilipschitz distortion, as in figure \ref{local spacing 
picture}, and it is an exercise in Euclidean geometry to see that the 
spacing is equal to the minimum of the spacings on either side. 

\realfig{local spacing picture}{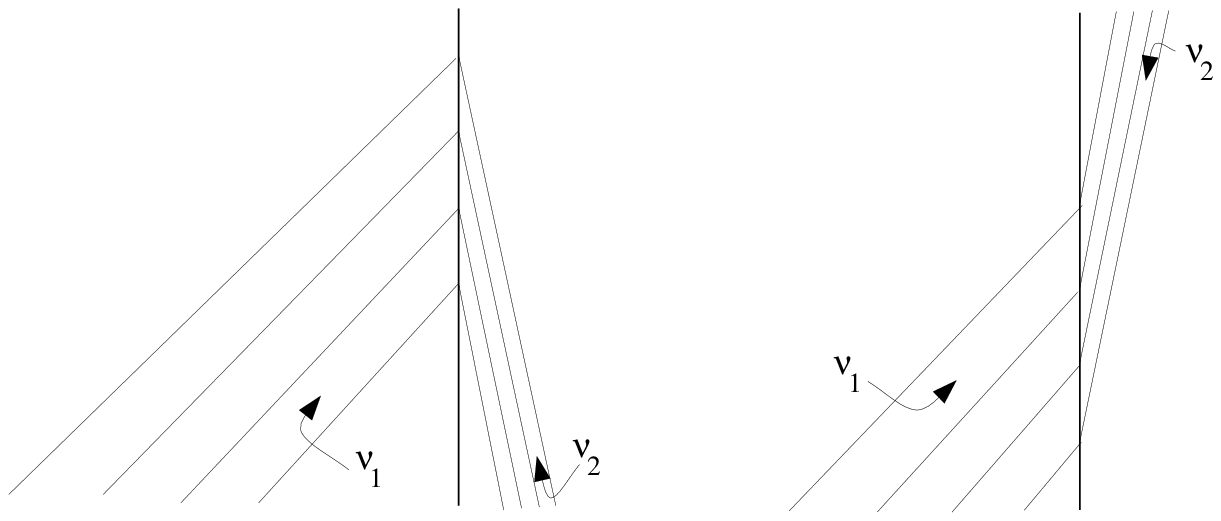}{}

Thus we have, since the area of each $Q_k$ is bounded, 
$$\nu_{Q_k,\sigma}(\beta_2\intersect Q_k) \ge {C\over 
\ell_\sigma(\beta^\sigma\intersect Q_k)}$$ 
for each $Q_k$, and for all of $P$ we obtain
$$
\nu_{P,\sigma}(\beta_2) \ge C'\min_k \nu_{Q_k,\sigma}(\beta_2\intersect Q_k)
$$
since again the spacing is not reduced by more than a constant factor at 
the interfaces between subsurfaces. 

Setting $\beta'=\beta_2$ and observing that 
$\ell_\sigma(\beta^\sigma\intersect P)\ge 
\max\limits_k \ell_\sigma(\beta^\sigma\intersect Q_k)$, the statement of 
lemma \ref{adjust on P} follows.
\end{proof}

Now we can apply the spacing lemma \ref{spacing estimate} and the fact 
that the area of $P$ is bounded in terms of its Euler characteristic
to conclude that 
$$
\lambda_{P,\sigma}(\beta) \le {C\over \nu^2_{P,\sigma}(\beta')} \le 
     C \ell^2_\sigma(\beta^\sigma\intersect P).
$$

To obtain the inequality in the other direction we observe that the 
minimal length in the homotopy class $\RR_P(\beta)$ is nearly attained by 
$\beta^\sigma\intersect P$. In fact, the shortest representative 
$\beta^P$ can be obtained by sliding the endpoints of 
$\beta^\sigma\intersect P$ around $\boundary P$ until an orthogonal arc 
is obtained. The sliding distance is uniformly bounded by lemma
\ref{bounded slide} which appears below.
Thus, $\ell_\sigma(\beta^P) \ge 
c\ell_\sigma(\beta^\sigma\intersect P)$ for uniform $c$, and we conclude
by the definition of extremal length that
$$
\lambda_{P,\sigma}(\beta) \ge c^2 {\ell_\sigma^2(\beta^\sigma)\over 
\Area(S,\sigma)}.
$$
Since $\Area(S,\sigma)$ is determined by the Euler characteristic of $S$, 
we are  done with lemma \ref{extremal is hyperbolic for thick}.
\end{proof}

\begin{lemma}{bounded slide}{}
Let $A_1$ and $A_2$ be annuli (possibly the same) of the $(\ep_0,\ep_1)$ 
collar decomposition of $(S,\sigma)$, and let $\alpha$,$\beta$ be geodesic arcs 
in $S-(A_1\union A_2)$ with endpoints on $\boundary A_1\union\boundary  
A_2$, such that $\alpha$ is homotopic to $\beta$ keeping the endpoints on 
the boundary. Then the homotopy can be taken to move the endpoints of 
$\alpha$ by no more than $c=c(\ep_0,\ep_1)$.
\end{lemma}

\begin{proof}{}
Lifting to $\til A_1,\til A_2,\til\alpha$ and $\til \beta$ in $\Hyp^2$, we 
see that the endpoints of $\til\alpha$ and $\til\beta$ must lie in the
arcs $a_1$  
and $a_2$ in figure \ref{bounded slide fig}. It is fairly easy to see 
that the extremal case (longest $a_i$) occurs when $\til A_1$ and $\til 
A_2$ are actually tangent, and each is the $r$-neighborhood of a 
geodesic, where $r=\cosh^{-1}(\ep_0/\ep_1)$. The length of $a_1$ and 
$a_2$ in this case gives the constant $c$.

\end{proof}

\realfig{bounded slide fig}{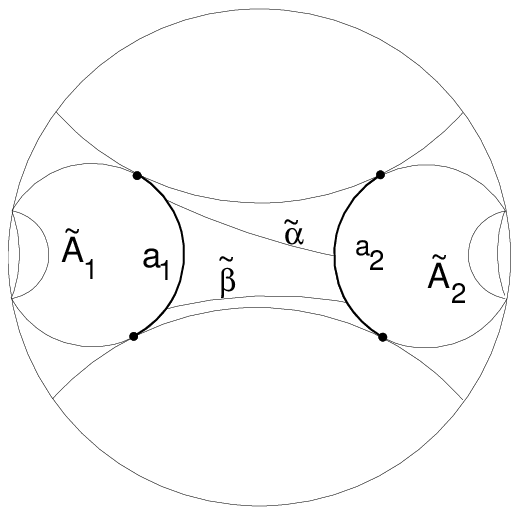}{}

\subsection{Extremal length on annuli.}
For annuli $A$ in $S$ we will define a quantity $\lambda_{A,\sigma}(\beta)$ which is 
not exactly extremal length, but will quantify the extremal length 
contribution of $A$ to the whole curve $\beta\in\CC(S,\boundary S)$.

In particular, if $A$ is an annulus of the $(\ep_0,\ep_1)$ collar 
decomposition of  
$(S,\sigma)$ with core curve $\gamma$, then for
$\beta\in\CC(S,\boundary S) $
let $i(\beta,A)$ or $i(\beta,\gamma)$ denote the (minimal) intersection 
number of  $\beta$ with $\gamma$ (note, this includes the case that
$A$ is peripheral, $\gamma\subset \boundary S$ and $\beta$ has
$i(\beta,\gamma)$ endpoints on $\gamma$).
The geodesic
$\beta^\sigma$ can meet $A$ either in a collection of $i(\beta,A)$  
properly embedded  
arcs, or  if $i(\beta,A)=0$, in a union of closed curves homotopic to
$\gamma$.
In either case let $n(\beta,A)$ denote the number of components of 
$\beta^\sigma\intersect A$. We define

\begin{equation}
\label{A length def}
\lambda_{A,\sigma}(\beta)  = 
\begin{cases}
i(\beta,A)^2 \left(m(A) + {t_{A,\sigma}(\beta)^2\over m(A)}\right)
& \text{if } i(\beta,A)>0, \\
n(\beta,A)^2/ m(A) & \text {otherwise.}
\end{cases}
\end{equation}

It is worth nothing that, if $t_{A,\sigma}(\beta)$ is an integer 
and a gluing together (without twisting) of the boundaries of $A$ takes 
$\beta$ to a closed curve on the torus, then (\ref{A length def}) gives 
the extremal length of the homotopy class of this curve. Compare also to 
lemma \ref{Kerckhoff-Royden for torus}.


%

\section{The main length estimate}
\label{main length estimate}
For the remainder of the paper we shall fix a choice of $\ep_0>\ep_1>0$ 
so that the following hold for any standard hyperbolic metric $\sigma$ on 
$S$.
\begin{enumerate}
\item $\ep_0$ is less than $\ep_\MM$, the Margulis constant of $\Hyp^2$.
\item Each annulus of the $(\ep_0,\ep_1)$ collar decomposition
has modulus at least 1. It is furthermore
contained in a slightly larger annulus, so that each internal boundary 
component bounds an annulus in the thick part of modulus $m_0$, where
$m_0$ is the modulus bound in lemma \ref{insulated surface}.
\item A simple geodesic in $(S,\sigma)$, if it intersects an annulus of the 
$(\ep_0,\ep_1)$ collar decomposition, does so either in arcs that connect the 
boundaries or it is the  geodesic core of the annulus.
\item Lemma \ref{internal and external twists} holds.
\end{enumerate}

Let $\sigma$ denote a fixed standard hyperbolic metric on $S$, and let 
$\QQ=\PP\union\AAA$ denote the $(\ep_0,\ep_1)$ collar decomposition.
of section \ref{thick-thin}. 

\begin{theorem}{main estimate}{}
Let $S$ be a surface of finite type with
a standard hyperbolic metric $\sigma$ in which
boundary lengths are at most $\ell_0$, and let $\QQ$ be the set of 
components of the $(\ep_0,\ep_1)$ collar decomposition.
Then, for any $\beta\in\CC(S,\boundary S)$ we have
$$
      \lambda_{S,\sigma}(\beta) \simmult 
                          \max_{Q\in\QQ} \lambda_{Q,\sigma}(\beta)
$$
where the constants for $\simmult$ depend only on $\ep_0,\ep_1,\ell_0$ 
and the topological type of $S$.
\end{theorem}

\begin{proof}{} 
Fix, for the remainder of the section, a class $\beta\in\CC(S,\boundary S)$.
Throughout the rest of the section we shall make the 
abbreviations
$\lambda_Q = \lambda_{Q,\sigma}(\beta)$, $t_A = 
t_{A,\sigma}(\beta)$, $m_A = m(A,\sigma)$, and $n(A)=n(\beta,A)$.
On each  annulus
$A\in\AAA$, let $\sigma_e$ denote the unique standard Euclidean metric
on $A$ which is conformally equivalent to $\sigma$, and agrees with
$\sigma$ on $\boundary A$.

\bfheading{The upper bound.} 
To bound $\lambda_{S,\sigma}(\beta)$ from above, we shall need to construct 
a metric on $S$, in the conformal class of $\sigma$, which is sufficiently 
similar to the extremal metric for $\beta$ that the spacing estimate in 
lemma \ref{spacing estimate} can be applied. We remark that it is of course 
known that an extremal metric exists, and in fact is given by a quadratic 
differential. Our construction, by contrast, is not extremal (in fact 
we make no explicit comparison to the extremal metric), but gives concrete 
estimates which are explicitly related to the hyperbolic geometry of 
$(S,\sigma)$. 

Let us make some new definitions.
Let $h_1:\QQ\to[0,\infty)$ be given by
\begin{equation}
\label{h1 def}
h_1(Q) = \left\{\begin{array}{ll}
			\sqrt{\lambda_Q} & \qquad Q\in\PP\\
			\sqrt{\lambda_Q/m_Q} & \qquad Q\in\AAA.
		  \end{array}\right.
\end{equation}	
This quantity has the property that the scaled metrics $h_1(P)\sigma$ on
$P\in\PP$ and $h_1(A)\sigma_e$ on $A\in\AAA$ have areas comparable to 
$\lambda_P$ and $\lambda_A$.

However, in order to scale the metric in a consistent way on the
entire surface we shall need to adjust $h_1$ to obtain a function $h$,
as follows:

\begin{lemma}{Scaling function}{}
There exists a positive function $h:\QQ\to\R_+$ such that the following 
hold.
\begin{enumerate}
\item $h(Q)\ge h_1(Q)$ for all $Q\in\QQ$,
\item $\max_Q h(Q) = \max_Q h_1(Q)$,
\end{enumerate}
and, whenever $A\in\AAA$ is adjacent to $P\in\PP$,
\begin{enumerate} \setcounter{enumi}{2}
\item $h(A)\le h(P)$.
\end{enumerate}
Furthermore, for any internal
annulus $A\in\AAA$ which is  adjacent to components $P_1,P_2\in\PP$:
\begin{enumerate} \setcounter{enumi}{3}
\item $\displaystyle  h(A)^2 \ge h(P_1)h(P_2) e^{-\pi m_A},$
\end{enumerate}
and if $A$ is peripheral and adjacent to $P\in\PP$,
\begin{enumerate} \setcounter{enumi}{4}
\item $\displaystyle h(A) \ge h(P)e^{-\pi m_A/2}.$
\end{enumerate}

\end{lemma}

\begin{proof}{}
Define $h_2(A) = h_1(A)$ for $A\in\AAA$. For $P\in\PP$, let
$A_1,\ldots,A_k\in\AAA$ be the annuli adjacent to $P$. Define
\begin{equation}
\label{h2 def thick}
h_2(P) =
\max\left\{{h_1(P)},{h_1(A_1)},\ldots,{h_1(A_k)}\right\}.
\end{equation}
This function satisfies conditions (1), (2) and (3) of the lemma, and note also 
that $h_2(P)>0$ for at least one $P$. 
Next, let us adjust $h_2$ to a function $h_n$ that satisfies the following
additional condition: for each internal $A\in\AAA$ adjacent to 
$P_1,P_2\in\PP$, we require that $h(P_1)h(P_2) > 0$, and 
\begin{equation}
\label{ratio condition}
\left|\log{h_n(P_1)\over h_n(P_2)}\right| \le {\pi m_A}.
\end{equation}
If the condition fails for $h_2$, then it must fail with some $A,P_1,P_2$
for which $h_2(P_1)>0$ and $h_2(P_2)<h_2(P_1)$. Define
\begin{equation}
\label{h3 def}
h_3(P_2) = h_2(P_1) e^{-\pi m_A},
\end{equation}
and $h_3=h_2$ elsewhere. Clearly $h_3\ge h_2$, and the maximal value
of $h_2$ hasn't been changed. We may repeat this process to obtain a
sequence of functions $h_1,h_2,h_3,\ldots$, which clearly terminates
after a finite number of steps in a function $h_n$ that satisfies 
conditions (1),(2),(3) and (\ref{ratio condition}). Note also that $h_n$ is 
positive on $\PP$, since if $h_n(P)=0$ after these adjustments, then 
$P$ must be adjacent only to peripheral annuli on which $h_1=0$, and in 
that case $\beta$ is empty (we are assuming that $S$ is connected).

For an internal annulus $A\in\AAA$ adjacent to $P_1,P_2$, if $h_n$ fails to
satisfy condition (4) we make a final adjustment:
\begin{equation}
\label{final h adjustment 1}
h(A) = \left(h_n(P_1)h_n(P_2)e^{-\pi m_A}\right)^{1/2}.
\end{equation} 
If $A$ is peripheral, adjacent to $P$ and fails condition (5), we adjust by
\begin{equation}
\label{final h adjustment 2}
h(A) = h_n(P)e^{-\pi m_A/2}
\end{equation} 
We let $h=h_n$ for other annuli and for $P\in\PP$. Note that $h$
still satisfies conditions (1)-(3). In particular, condition (3) is
satisfied by virtue of the inequality (\ref{ratio condition}). This
concludes the proof of lemma \ref{Scaling function}.
\end{proof}

The next step is to define  a new metric $\sigma'$ on $S$, conformally
equivalent to $\sigma$.
For $P\in\PP$ we have
simply
\begin{equation}
\label{metric def thick}
\sigma'|_P = h(P)\sigma.
\end{equation}
For an internal annulus $A$, let $P_1$ and $P_2$ be
the components of $\PP$ which adjoin  $A$. For each $P_i$, if
$h(P_i)>h(A)$ we must interpolate between the natural scalings for
$P_i$ and $A$. Let $m=m_A$ and let
\begin{equation}
\mu_i = {1\over 2\pi}\log(h(P_i)/h(A)).
\end{equation}
This is the
modulus of an annular region in the plane between two
concentric circles whose radii are in the ratio $h(P_i)/h(A)$. We
would like $\sigma'$ to be precisely such a metric on subannuli of $A$
of appropriate modulus. In addition we want
to reserve a sufficiently large (modulus at least $m/2$) middle section of
$A$  for other purposes. It is for these reasons that condition (4)
was imposed in lemma \ref{Scaling function}.

Condition (4) implies that  $\mu_1+\mu_2\le m/2$. 
Thus, let $B_1$ and
$B_2$ be annulus neighorhoods of the boundary components of $A$, whose
radii in the Euclidean metric $\sigma_e$ are $\ep_0\mu_1,\ep_0\mu_2$. Let $C$
denote the closure of the complementary annulus $A-B_1-B_2$; the
modulus of $C$ is at least $m/2$. 

On $B_1$ and $B_2$ we can make 
$\sigma'$ the metric conformally equivalent to $\sigma$ such that
$(B_i,\sigma')$ is a planar annulus as above, with boundary lengths
$\ep_0h(P_i)$ (for the boundary between $A$ and $P_i$) and
$\ep_0h(A)$ on the inner boundary. On the remaining annulus $C$ between
$B_1$ and $B_2$ 
let $\sigma'= h(A)\sigma_e$. Note that this makes $\sigma'$ a
continous metric on $S$. (See figure \ref{metric in pieces}.)

\realfig{metric in pieces}{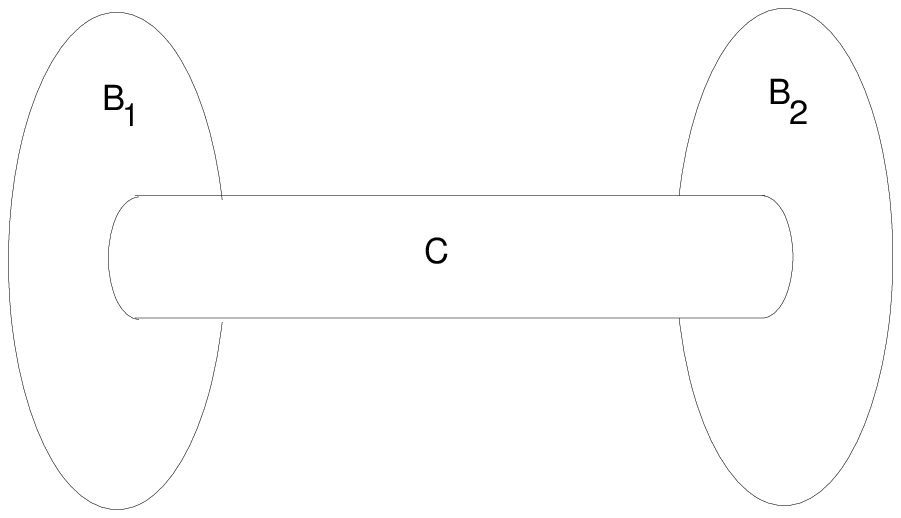}{The metric
$\sigma'$ on the annulus $A$ decomposed into $B_1,C$ and $B_2$.}

Similarly if $A$ is a peripheral annulus adjacent to $P$, we can apply 
condition (5) of lemma \ref{Scaling function} implies that 
$\mu = {1\over 2\pi}\log(h(P)/h(A)) \le m/2$. Thus we place a planar 
annulus metric on the annulus neighborhood $B$ of $P\intersect A$, of 
$\sigma_e$-radius $\ep_0 \mu$. On $A-B$ we again let $\sigma' = 
h(A)\sigma_e$. 

\bfheading{Finding a good representative.}
Having defined $\sigma'$, we shall now construct a
representative of $\beta$ whose spacing in $\sigma'$ is as uniform as
possible. 

\begin{lemma}{uniformly spaced representative}{}
There is a representative $\beta_2$ of $\beta$ whose spacing in the
metric $\sigma'$ is bounded by
$$
	\nu_{S,\sigma'}(\beta_2) \ge c.
$$
The constant $c$ depends only on $\ep_0,\ep_1,\ell_0$ and $\chi(S)$.
\end{lemma}

\begin{proof}{}
We begin by applying lemma \ref{adjust on P} to deform $\beta\intersect 
P$ (if non-empty) on each $P\in\PP$ to a curve system
$\beta_1$ whose intersection with each boundary component of $P$ 
is evenly spaced, and whose spacing on the interior satisfies
$\nu_{P,\sigma} \ge c/\ell_\sigma(\beta^\sigma\intersect P)$.
(The portions of $\beta$ in the annuli $A\in\AAA$ can be moved along with 
the homotopy on $P\in\PP$, and we will examine these separately).
By lemma \ref{extremal is hyperbolic for thick}, this implies
$$\nu_{P,\sigma}(\beta_1) \ge c'/\sqrt{\lambda_P} = c'/h_1(P)$$
and therefore
$$
    \nu_{P,\sigma'}(\beta_1) \ge c' h(P)/h_1(P) \ge c'.
$$

Now consider the case of an internal annulus $A\in\AAA$, partitioned as 
before into 
boundary annuli $B_1,B_2$ and a central annulus $C$. 
Suppose first that $i(\beta,A)>0$.
Deform $\beta_1\intersect A$, fixing its intersection with
$\boundary A$, so that each resulting arc intersects $(B_i,\sigma')$
in a radial segment, and $(C,\sigma')$ in a $\sigma'$-geodesic arc that twists
several times around $C$. Call this new curve $\beta_2$. We 
have
\begin{equation}
t_{C,\sigma'}(\beta_2\intersect C) = t_A + \delta
\end{equation}
where the error satisfies $|\delta|\le 1$, since
by lemma \ref{adjust on P} the endpoints have moved at most halfway 
around the boundary in the homotopy from $\beta$ to $\beta_1$. 
It follows that the length of
$\beta_2\intersect C$ is, if $i(\beta,A)>0$, 
\begin{eqnarray}
\label{length for pos i}
\ell_{\sigma'}(\beta_2\intersect C) &=& n(A)h(A)\ep_0\sqrt{m(C)^2 +
(t_A+\delta)^2 } \nonumber \\ 
&\le & c\ep_0h(A)\sqrt{m_A\lambda_A}
\end{eqnarray}
Here we are using the formula (\ref{A length def})
for $\lambda_A=\lambda_{A,\sigma}(\beta)$.
Note that the error $\delta$ can be replaced by the multiplicative factor
$c$ because $m_A\ge 1$ (by choice of $\ep_0,\ep_1$).

If $i(\beta,A)=0$ and $n(A)>0$ we can make $\beta_2\intersect A$ a 
disjoint union 
of closed $\sigma'$-geodesics uniformly spaced in $C$. Then we have:
\begin{eqnarray}
\label{length for zero i}
\ell_{\sigma'}(\beta_2\intersect C) &=& n(A)h(A)\ep_0 \nonumber\\
				&=& \ep_0 h(A) \sqrt{m_A\lambda_A}
\end{eqnarray}
again using (\ref{A length def}).
Now, since $\beta_2$ is uniformly spaced in $C$ with respect to 
the standard Euclidean metric $\sigma'$, its spacing is equal to the area 
divided by its length. 
Applying either (\ref{length for pos i}) or (\ref{length 
for zero i}), and the fact that $m(C)\ge m_A/2$, we obtain:
\begin{eqnarray}
\nu_{C,\sigma'}(\beta_2\intersect C) &=& {\Area(C,\sigma')\over
				\ell_{\sigma'}(\beta_2\intersect C)}
								\nonumber\\ 
			&\ge& c{\ep_0^2h(A)^2m_A/2\over
			\ep_0 h(A) \sqrt{m_A\lambda_A}} \nonumber\\
			&=& c' {h(A)\over h_1(A)}.
\end{eqnarray}
The spacing in the boundary annuli $B_1,B_2$ is no smaller than this,
and $h(A)\ge h_1(A)$ by lemma \ref{Scaling function}, 
so we conclude 
\begin{equation}
\nu_A \ge c'.
\label{Spacing bound in A}
\end{equation}

The case when $A$ is peripheral is dealt with in a similar way, and 
(\ref{Spacing bound in A}) still holds. In fact it is slightly easier 
since $t_A=0$ in this case. 

Thus we have a curve $\beta_2$ with a uniform lower bound on $\nu_{Q,\sigma'}$ 
for all $Q\in\QQ$. We can now make the same observation as in the proof 
of lemma \ref{adjust on P} to conclude that 
$$\nu_{S,\sigma'}(\beta_2) \ge c\min_{Q\in\QQ} \nu_{Q,\sigma'}(\beta_2).$$
Namely, since on a uniformly sized neighborhood of an interface between 
two components the situation is uniformly bilipschitz equivalent to the 
Euclidean situation depicted in figure \ref{local spacing picture}, the 
spacing of the  combined curve is no less than a constant times 
the minimum of the spacings on either side.  
This concludes the proof of lemma \ref{uniformly spaced representative}.
\end{proof}

We can now apply lemma \ref{spacing estimate} to obtain:
\begin{equation}
\label{area bounds lambda}
\lambda_{S,\sigma}(\beta) \le c\Area(S,\sigma') = c\sum_{Q\in\QQ} 
\Area(Q,\sigma').
\end{equation}
The desired upper bound on $\lambda_{S,\sigma}$ will be obtained if we 
show that 
\begin{equation}
\label{lambda bounds area}
\max_Q \Area(Q,\sigma') \le c\max_Q \lambda_Q
\end{equation}
for a uniform constant $c$, since the number of terms in the sum
(\ref{area bounds lambda}) is bounded in terms of the topology of $S$.

For any $P\in\PP$, recall that 
$\Area(P,\sigma') \le ch^2(P)$ for a constant $c$. An internal annulus
$A\in\AAA$ is divided into three subannuli $B_1,B_2$ and $C$ whose areas
are $\Area(C,\sigma') \le \ep_0 h^2(A) m_A$, and 
$\Area(B_i,\sigma') = \ep_0^2 (h^2(P_i)-h^2(A))/4\pi$, where $P_1,P_2$
are the neighboring components of $\PP$. A peripheral annulus is treated 
similarly. Thus, it suffices to bound
the quantities $h^2(P)$ and $h^2(A)m_A$.

Note first that, by definition (\ref{h1 def}), $h^2_1(P) \le
\lambda_Q$, since $m_A\ge 1$ for all $A\in\AAA$. Thus we have
\begin{equation}
\label{max h1}
\max_{Q\in\QQ} h_1^2(Q) \le \max_Q\lambda_Q.
\end{equation}
It follows from this and part (2) of lemma \ref{Scaling function} that,
for any $P\in\PP$,
\begin{equation}
\label{h(P) bound}
h^2(P) \le \max_Q\lambda_Q.
\end{equation}
Now for $A\in\AAA$, suppose first that $h(A)=h_1(A)$. In this case
\begin{align}
h^2(A)m_A &= h_1^2(A) m_A \notag\\
	& = \lambda_A \le \max_Q\lambda_Q
\label{h(A) bound 1}
\end{align}
If $h(A)>h_1(A)$ then, if $A$ is internal we use (\ref{final h 
adjustment 1}) in the proof of lemma \ref{Scaling function} to obtain
\begin{align}
h^2(A)m_A 	& = h_n(P_1)h_n(P_2) m_A e^{-\pi m_A} \notag\\
		& \le \max_Q h_1^2(Q)  m_A e^{-\pi m_A}  \notag\\
		& \le {1\over\pi  e}\max_Q \lambda_Q 
\label{h(A) bound 2}
\end{align}
where the last line follows from (\ref{max h1}) and the fact that
$xe^{-\pi x} \le 1/\pi e$ for all $x$. In case $A$ is peripheral we apply 
(\ref{final h adjustment 2}), to similar effect.
The bounds (\ref{h(P) bound}-\ref{h(A) bound 2})
together imply 
(\ref{lambda bounds area}),  so that with (\ref{area
bounds lambda}) we have the desired upper bound
$$
\lambda_{S,\sigma}(\beta) \le c\max_{Q\in\QQ} \lambda_{Q,\sigma}(\beta).
$$

\bfheading{The lower bound.}
For each $Q\in\QQ$, we shall find a metric $\sigma''$ conformal to 
$\sigma$ such that $\ell^2_{\sigma''}([\beta])/\Area(S,\sigma'')$ is 
approximately $\lambda_Q$. This will give the desired lower bound. 

If $Q=P\in\PP$, let $\sigma'' = \sigma$. Note that
$\ell_\sigma([\beta])\ge\ell_\sigma(\beta^\sigma\intersect P) \simmult 
\sqrt{\lambda_Q}$, and that $\Area(S,\sigma)$ is a constant, 
$2\pi|\chi(S)|$. Thus: 
\begin{equation}
\lambda_{S,\sigma}(\beta)\ge {\ell_\sigma^2(\beta)\over\Area(S,\sigma)} \gesim
\lambda_Q
\end{equation}

If $Q=A\in\AAA$, define the metric $\sigma^A$ to be 
\begin{equation}
\sigma^A = \left\{
\begin{array}{ll}
	\sigma & \hbox{in $S-A$}  \\
	\sigma_e & \hbox{in $A$}
\end{array}	\right.
\end{equation}
where $\sigma_e$ is the Euclidean metric defined as in the beginning of 
this section. 
 
If $i(\beta,A) = 0$ but $n(A)>0$, there are $n(A)$ copies in $\beta$ of 
the core curve $\gamma$ of $A$. It is easy to see that 
$\ell_{\sigma^A}(\gamma)  = \ep_0$ in this case, since $\sigma_A$ is 
non-positively curved so that the $\sigma_A$-geodesic cores of $A$ are 
the shortest representatives of $\gamma$. Thus, 
$$
\lambda_{S,\sigma}(\beta) \ge {n^2(A)\ep_0^2 \over \ep^2_0 m_A + 
2\pi|\chi(S)|} \gesim \lambda_A.
$$
where for the last step we use the fact that $m_A\ge 1$ to
replace an additive constant with a multiplicative one. 

Now consider the case where $i(\beta,A)>0$.
Let $\beta^A$ be a representative of $\beta$ of minimal $\sigma^A$-length
(in fact $\beta^A$ is 
unique, but we do not need to use this), and note that each arc of 
$\beta^A\intersect (S-A)$ 
is geodesic, and homotopic with endpoints on $\boundary A$ to a 
corresponding arc of $\beta^\sigma\intersect (S-A)$ (this is most clearly seen 
in the universal cover). Thus by lemma \ref{bounded slide}, the homotopy 
between $\beta^\sigma$ and $\beta^A$ only moves the intersection with 
$\boundary A$ a bounded amount $c$. In particular 
we may conclude that 
$$|t_{A,\sigma}(\beta^A) - t_A| \le 2c/\ep_0,$$
and therefore
\begin{eqnarray}
	\ell_{\sigma^A}(\beta^A\intersect A) & \ge  & 
	n(A)\ep_0\sqrt{m^2(A)+(t_A-2c/\ep_0)^2} \nonumber  \\
	 & \gesim & n(A)\ep_0 \sqrt {m^2(A)+t^2_A} \nonumber \\
     & \gesim	 & \sqrt{\lambda_A m_A}. 
\end{eqnarray}
(Note that we have again used the fact that $m_A\ge 1$.)

We conclude that 
\begin{equation}
\lambda_{S,\sigma}(\beta) \gesim {\lambda_A m_A \over \ep_0 m_A + 
2\pi|\chi(S)|} \gesim \lambda_A.
\end{equation}
This completes the proof of the lower bound and therefore of theorem
\ref{main estimate}
\end{proof}

We can generalize theorem \ref{main estimate} slightly as follows. 
Suppose that $A_1,\ldots,A_k$ are a subset of the annuli of the $(\ep_0,\ep_1)$ 
collar decomposition of $(S,\sigma)$ and $P_1,\ldots,P_m$ are the 
components of $S-\union_i A_i$. We call $\QQ = 
\{A_1,\ldots,A_k,P_1,\ldots,P_m\}$ a {\em partial} $(\ep_0,\ep_1)$ 
decomposition. We can then state: 

\begin{corollary}{generalized estimate}{}
If $\QQ$ is a partial $(\ep_0,\ep_1)$ collar decomposition with 
$\ep_0,\ep_1$ chosen as above, then for any $\beta\in\CC(S,\boundary S)$ 
$$
\lambda_{S,\sigma} \simmult \max_{Q\in\QQ} \lambda_{Q,\sigma}(\beta).
$$
\end{corollary}

\begin{proof}{}
Let $\bar\QQ=\bar\PP\union\bar\AAA$ be the full $(\ep_0,\ep_1)$ collar 
decomposition.  Let $\bar\QQ_i = \bar\PP_i\union \bar\AAA_i$ be the 
subset of $\bar \QQ$  consisting of components that lie  in $P_i$. 
Let $\hat\sigma_i$ be the standard hyperbolic metric on $P_i$ which is 
conformally equivalent to the restriction of $\sigma$. 
We first show that, for any $\bar Q\in\bar \QQ_i$,  
\begin{equation}
\label{subsurface ok}
\lambda_{\bar Q,\sigma}(\beta) \simmult 
\lambda_{\bar Q,\hat\sigma}(\RR_{P_i}(\beta)).
\end{equation}
For a hyperbolic component $\bar Q \in \bar\PP_i$, these two quantities are  
the same by definition (\ref{P length def}), since only the conformal 
class of the metric matters for extremal length. 

For an annulus $\bar Q=\bar A\in\bar\AAA_i$, the definition (\ref{A
length def})  
depends in  particular on the measure of twisting $t_{\bar A,\sigma}$, which 
depends on the actual geometry. However, lemma 
\ref{twist in subsurface} assures us that the difference between twisting 
in $\sigma$ and $\hat\sigma$ is no more than a constant. 
It follows that (\ref{subsurface ok}) holds for annuli as well, with 
uniform constants that depend only on $\ep_0,\ep_1$ and the topology of 
$S$. 

The corollary now follows by applying theorem \ref{main estimate} 
separately to each $P_i$.

\end{proof}

\section{The product region theorem}
\label{product structure}

In this section $S$ denotes an oriented finite genus surface without 
boundary, but possibly with finitely many punctures.

Let $\gamma$ be an element of $\CC(S)$ all of whose components 
$\gamma_1,\ldots,\gamma_k$ are homotopically distinct, and fix an 
orientation on each component of $\gamma$.
For any $\ep>0$, denote by
$Thin_\ep(S,\gamma)$ the subset of $\TT(S)$ consisting of all $[\sigma]$ 
such that $\ell_{\sigma}(\gamma_i)\le \ep$. 

Let $\gamma_{k+1},\ldots,\gamma_n$ denote a completion of $\gamma$ to
an oriented pair-of-pants decomposition $\hat\gamma$. Fix also a curve
system $\mu$ composed of seam arcs, giving a marking of $S$ as in
\S\ref{twisting section}. This gives rise to length and twist
coordinates on $\TT(S)$, denoted in \S\ref{twisting section} by the
map $F:\TT(S)\to \R_+^n\times\R^n$. Let us write the coordinates of
$F$ as $F_i$, where $F_1,\ldots,F_n$ are length coordinates, and
$F_{n+1},\ldots,F_{2n}$ are twist coordinates.
Separating out the coordinates of the first $k$ curves
from the rest decomposes $\TT(S)$ as a product. 

Let $S_\gamma$ denote the punctured (or ``noded'') surface obtained
from $S$ by deleting $\gamma_1,\ldots,\gamma_k$ 
and replacing each curve by a pair of punctures. Let $\TT(S_\gamma)$   
be the space of (analytically finite) conformal structures on $S_\gamma$. 
Let $\Hyp$ denote the upper half-plane $\{(x,y):y>0\}$ endowed with the 
hyperbolic metric $d_\Hyp$. Let $X_\gamma$ denote
the product space
$X_\gamma = \TT(S_\gamma)\times\Hyp_1\times\cdots\times\Hyp_k$ where each 
$\Hyp_i$ is a copy of $\Hyp$, endowed with the sup metric
$d_X= \max(d_{\TT(S_\gamma)},d_{\Hyp_1},\ldots,d_{\Hyp_k})$.

The marking by $\hat\gamma$ and $\mu$ induces a marking on each
component of $S_\gamma$, and a map
$G:\TT(S_\gamma)\to \R_+^{n-k}\times\R^{n-k}$, given by
$G(\rho) =
(\ell_\rho(\gamma_{k+1}),\ldots,\ell_\rho(\gamma_{n}),
s_{k+1}(\rho),\ldots,s_n(\rho))$ for any complete hyperbolic metric
$\rho$ on $S_\gamma$.

We can therefore define a map
\begin{equation}
\Pi: \TT(S) \to \TT(S_\gamma)\times 
\Hyp_1\times\cdots\times\Hyp_k
\label{Pi def}
\end{equation}
as follows.  Let $\Pi_0$ denote the component of $\Pi$ mapping to
$\TT(S_\gamma)$, and for $i>0$ let $\Pi_i$ denote the component
mapping to $\Hyp_i$.  We define $\Pi_0$ via
\begin{equation}
\label{Pi0 def}
\Pi_0(\sigma) =
G^{-1}(F_{k+1}(\sigma),\ldots,F_{n}(\sigma), 
F_{n+k+1}(\sigma),\ldots,F_{2n}(\sigma)).
\end{equation}
In other words,
we build a punctured surface by forgetting the length and twist
coordinates of the first $k$ curves, and using the remaining
coordinates to build new pairs of pants. 
We then use the forgotten coordinates for the rest of $\Pi$, by defining
\begin{equation}
\label{Pii def}
\Pi_i(\sigma) = (s_i(\sigma),1/\ell_\sigma(\gamma_i))
\end{equation}
for $i=1,\ldots,k$.

Let $\ep_0$ and $\ep_1$ be chosen as in the previous section. Our main 
theorem is the following.

\begin{theorem}{Product structure}{}
Let $\gamma\in\CC(S)$ be an oriented curve system as above, extended
to a marking $\hat\gamma,\mu$ of $S$. Given
$\ep\in(0,\ep_1)$, there is a 
contant $a_0$ depending on $\ep_0,\ep_1$ and the topological 
type of $S$, such that
the homeomorphism $\Pi:\TT(S) \to X_\gamma$ 
defined above, when restricted to $Thin_\ep(S,\gamma)$, 
satisfies
$$
|d_{\TT(S)}(\sigma,\tau) - d_X(\Pi(\sigma),\Pi(\tau))| \le a_0
$$
for any $[\sigma],[\tau] \in Thin_\ep(S,\gamma)$.
\end{theorem}

\begin{proof}{}
For a given $[\sigma]\in Thin_\ep(S,\gamma)$, let $\sigma$ denote its
hyperbolic representative. 
Let $A_1,\ldots,A_k$ be the annuli of the $(\ep_0,\ep_1)$ collar 
decomposition which 
correspond to the curves $\gamma_1,\ldots,\gamma_k$. The complement 
$S-\union A_i$ is also naturally identified with $S_\gamma$.
%
Let $P_1,\ldots,P_m$ denote the connected components of $S_\gamma$,
and let $\QQ =\{A_j\}\union\{P_i\}$. Note that for another point
$[\tau]\in Thin_\ep(S,\gamma)$ we obtain a similar decomposition,
which we also call $\QQ$. 

Our strategy will be to estimate distances using Kerckhoff's theorem, so 
we shall need to be able to control extremal length ratios of the type
$\lambda_{Q,\sigma}(\alpha)/\lambda_{Q,\tau}(\alpha)$. 
Note that, for $Q$ a subsurface of $S$, if $\lambda_{Q,\sigma}(\alpha)
= 0$ then $\lambda_{Q,\tau}(\alpha) = 0$ as well. 
Thus we adopt the convention that, if $\lambda_{Q,\sigma}(\alpha) = 0$
then the expression
$\lambda_{Q,\sigma}(\alpha)/\lambda_{Q,\tau}(\alpha)$ is defined to be 1.

The first step is  the following: 

\begin{lemma}{closed curves give ratio}{}
Let $Q$ be any hyperbolic surface of finite type with boundary, and
$\lambda_1>\lambda_0>0$ two given constants. 
Let $[\sigma],[\tau]\in \TT(Q)$ be two conformal structures such that 
$\lambda_{Q,\sigma}(\gamma) \in [\lambda_0,\lambda_1]$
for any boundary component $\gamma$, and similarly for $\tau$. Then
\begin{equation}
\sup_{\alpha\in\CC_0(Q)} 
{\lambda_{Q,\sigma}(\alpha)\over\lambda_{Q,\tau}(\alpha)}
\simmult
\sup_{\alpha\in\CC(Q,\boundary Q)} 
{\lambda_{Q,\sigma}(\alpha)\over\lambda_{Q,\tau}(\alpha)}
\label{sup closed equals sup relative}
\end{equation}
where the constants for ``$\simmult$" depend only on
$\chi(S),\lambda_0$ and $\lambda_1$.
\end{lemma}

\begin{proof}{}
The left side is clearly no larger than the right, since
$\CC_0(Q)\subset\CC(Q,\boundary Q)$. 

The supremum on the right side is equal to $\exp 2d_{\TT(Q)}(\sigma,\tau)$,
by Kerckhoff's theorem \ref{Kerckhoff distance thm},
where $d_{\TT(Q)}$ is the distance in the Teichm\"uller space of 
$Q$. In particular, this distance is  
symmetric, so that the supremum of the reciprocal ratio is the same 
quantity. Choose $\beta\in\CC(Q,\boundary Q)$ such that 
$\lambda_{Q,\tau}(\beta)/\lambda_{Q,\sigma}(\beta)$ is nearly equal to 
the supremum on the right side of (\ref{sup closed equals sup relative}). 

By lemma \ref{opposite intersection number lemma}, whose proof 
appears below,
there exists $\alpha\in\CC_0(S)$ such that, for a uniform
$c>0$, 
\begin{equation}
\label{opposite intersection number}
i(\alpha,\beta)^2 \ge c
\lambda_{Q,\tau}(\alpha)\lambda_{Q,\tau}(\beta).
\end{equation}
(One should think of $\alpha$ as approximately ``orthogonal'' to 
$\beta$.)

On the other hand, the following inequality is elementary for any two
curves $\alpha,\beta\in\CC(Q)$ and $[\sigma]\in\TT(Q)$ (see e.g. 
\cite[Lemma 5.1]{minsky:slowmaps}).
\begin{equation}
\label{easy intersection number}
\lambda_{Q,\sigma}(\alpha)\lambda_{Q,\sigma}(\beta) \ge
i(\alpha,\beta)^2.
\end{equation}

Combining (\ref{opposite intersection number}) with(\ref{easy 
intersection number}), we obtain
$$
{\lambda_{Q,\sigma}(\alpha)\over\lambda_{Q,\tau}(\alpha)} \ge
c { \lambda_{Q,\tau}(\beta) \over
    \lambda_{Q,\sigma}(\beta)}
$$
which gives the remaining direction of lemma 
\ref{closed curves give ratio}. 
\end{proof}

To complete lemma \ref{closed curves give ratio} we need to prove the 
following: 

\begin{lemma}{opposite intersection number lemma}
Let $Q$ be a hyperbolic surface of finite type, with a conformal 
structure $[\tau]$. Suppose that for $0<\lambda_0<\lambda_1$ we have 
$\lambda_{Q,\tau}(\gamma) \in [\lambda_0,\lambda_1]$ for any boundary 
component $\gamma$ of $Q$. Then for any 
$\beta\in\CC(Q,\boundary Q)$ there exists $\alpha\in\CC_0(Q)$ for which 
$$
i(\alpha,\beta)^2 \ge c\lambda_{Q,\tau}(\alpha) \lambda_{Q,\tau}(\beta)
$$
where $c$ depends only on $\lambda_0,\lambda_1$ and $\chi(Q)$. 
\end{lemma}

\begin{proof}{}
Represent $[\tau]$ by a standard hyperbolic metric $\tau$, and
let $T_1,\ldots, T_j$ be the $(\ep_0,\ep_1)$ collar 
decomposition of $Q$, with respect to $\tau$. Note that $\ep_1$ may be
chosen (depending on $\lambda_0$) so that no annuli in the
decomposition are peripheral in $Q$.
Applying theorem \ref{main estimate}, there exists one component, say $T_1$, 
such that $\lambda_{Q,\tau}(\beta) \simmult \lambda_{T_1,\tau}(\beta)$.

Suppose first that $T_1$ is a non-annular component. We claim that there 
exist constants $c_0,\lambda_2$ such that there is an $\alpha\in\CC_0(T_1)$ with 
$\lambda_{T_1,\tau}(\alpha)\le \lambda_2$, and $i(\alpha,\beta)^2 \ge 
c_0\lambda_{T_1,\tau}(\beta)$. We can see this as follows: 
Since $T_1$ has no thin parts we can cut it along moderate-length curves 
to obtain bounded-diameter pieces. More precisely, note that the 
components of $\boundary T_1$ are either boundaries of thin parts of 
length $\ep_0$, or boundaries of $Q$, which have $\tau $-length in some 
interval $[\ell_0,\ell_1]$ (depending on $\lambda_0,\lambda_1$). 
We can find a pair-of-pants decomposition for $T_1$ of minimal total 
boundary lengths, and then add sufficiently many simple closed geodesics, as 
short as possible, to cut $T_1$ into disks and boundary-parallel annuli. 
The result is a collection of  non-peripheral
simple closed geodesics $\alpha_1,\ldots,\alpha_q$ in $T_1$ ($q$ 
dependent only on  $\chi(T_1)$) such that 
$\ell_\tau(\alpha_i)\le \ell_2$ (depending on 
$\ep_0,\lambda_0,\lambda_1,\chi(T_1)$)
and which cut $T_1$ up into components that 
are either disks, or annuli which share a boundary with $T_1$, or
annular neighborhoods of cusps. In the first two cases the diameters
of the components are uniformly bounded. In the third case we may
remove a neighborhood of the cusp avoided by $\beta$ (since $\beta$ is
simple) and so that what remains has bounded diameter.
Furthermore note that any component of the intersection of the geodesic 
$\beta^\tau$ with one of these regions has length at most some $\delta_0$, 
since $\beta^\tau$ meets $\boundary Q$ orthogonally, and meets the other 
components of $\boundary T_1$ nearly orthogonally (compare lemma 
\ref{bounded slide}).

It follows, adding over the pieces of $\beta^\tau$ in these components, 
that
$\sum_j \delta_0i(\beta,\alpha_j) \ge \ell_\tau(\beta^\tau\intersect T_1)$.
By lemma \ref{extremal is hyperbolic for thick},
$\ell^2_\tau(\beta^\tau\intersect T_1) \simmult
\lambda_{T_1,\tau}(\beta)$. Thus, there must be some 
$\alpha = \alpha_j$ such that 
$$i(\alpha,\beta)^2 \ge C \lambda_{T_1,\tau}(\beta).$$
We also know that there is an upper bound
$$\lambda_{Q,\tau}(\alpha) \le \lambda_3$$
for a uniform $\lambda_3$, since there is an upper bound on the 
hyperbolic length of $\alpha$ by construction.
The statement of the lemma follows in this case. 

Now let $T_1=A$ be an annular component with modulus $m=m(A,\tau)$, and 
consider first the case  
when  $i(\beta,A)>0$, and $t_1 = t_{A,\tau}(\beta) \ne 0$.
Let $t_2 = -m^2/t_1$ -- we choose this number because two geodesics with 
twists $t_1,t_2$ in a standard Euclidean annulus of modulus $m$ are 
orthogonal when $t_1t_2 = -m^2$. 
Lemma \ref{closed curve from annulus segment} guarantees the existence of  
an element $\alpha\in\CC_0(Q)$  with $|t_{A,\tau}(\alpha) - t_2| \le 
c_1$ for a uniform $c_1$. 

Now using 
the definition of $\lambda_{A,\tau}$ in (\ref{A length def}), 
we can see that
\begin{equation}
\begin{split}
\lambda_{A,\tau}(\alpha)\lambda_{A,\tau}(\beta)
                        &\simmult n(\alpha,A)n(\beta,A)(m+t_1^2/m)(m+t_2^2/m) \\
& = n(\alpha,A)n(\beta,A)(t_1-t_2)^2. 
\label{length product beats twist difference}
\end{split}
\end{equation}
where the additive error $c_1$ becomes a multiplicative error since $m\ge 
1$. 

On the other hand it is easy to see that, if two 
properly embedded arcs
$\alpha_1,\alpha_2$ in $A$ have twisting numbers $t_i = 
t_{A,\tau}(\alpha_i)$, then their intersection number is estimated by
\begin{equation}
|t_1-t_2|-1\le i(\alpha_1,\alpha_2) \le |t_1-t_2|+1.
\label{intersection is twist difference}
\end{equation}
Combining (\ref{intersection is twist difference}) with (\ref{length 
product beats twist difference}) we obtain
\begin{equation}
\lambda_{A,\tau}(\alpha)\lambda_{A,\tau}(\beta) 
\simge i(\alpha,\beta)^2
\end{equation}
(again an additive error has become multiplicative).
By lemma \ref{closed curve from annulus segment}, $\alpha$ may be
chosen to have bounded length outside $A$, and thus by theorem
\ref{main estimate}, $\lambda_{A,\tau}(\alpha) \simmult 
\lambda_{Q,\tau}(\alpha)$. 
Since also $\lambda_{A,\tau}(\beta) \simmult \lambda_{Q,\tau}(\beta)$ by 
choice of $A$,
the statement of lemma \ref{opposite intersection number lemma} 
follows.

We have left the case where $t_1 = 0$, and the case where $i(\beta,A)=0$ 
but $n(\beta,A) >0$. In the former case, We can pick $\alpha$ to be the 
core of $A$ (essentially this is $t_2 = \infty$), and note that 
$\lambda_{Q,\sigma}(\alpha) \simmult 1/m$, so that the lemma holds automatically. 
In the latter case we may choose $\alpha$ with $t_2 = 0$, and the roles 
of $\alpha$ and $\beta$ are reversed. 
\end{proof}

We now resume the proof of theorem \ref{Product structure}.
Let us establish the following
for $[\sigma],[\tau]\in Thin_\ep(S,\gamma)$:
\begin{equation}
\label{Reduce to relative subsurface}
\sup_{\alpha\in\CC(S)} {\lambda_{S,\sigma}(\alpha)\over
\lambda_{S,\tau}(\alpha)} 
\simmult
\max_{Q\in\QQ}  \sup_{\alpha\in\CC(S)} {
		\lambda_{Q,\sigma}(\alpha) \over
		\lambda_{Q,\tau}(\alpha)  }.
\end{equation}
One direction is a direct application of the main extremal length 
estimate. Note that $\QQ$ is what we called a ``partial
$(\ep_0,\ep_1)$ collar decomposition'' in the previous section. 
Thus by corollary 
\ref{generalized estimate} of theorem \ref{main estimate}, we have
\begin{equation}
\label{Main length ratio}
\lambda_{S,\sigma}( \alpha) \simmult \max_{Q\in\QQ} 
\lambda_{Q,\sigma}(\alpha)
\end{equation}
for each $\alpha\in\CC(S)$. Letting $Q_\alpha$ maximize  
$\lambda_{Q,\sigma}(\alpha)$ over all $Q\in\QQ$, 
we apply (\ref{Main length ratio}) to obtain
\begin{eqnarray}
\sup_{\alpha\in\CC(S)}{\lambda_{S,\sigma}(\alpha)\over\lambda_{S,\tau}(\alpha)}
		& \simmult &  
			\sup_{\alpha\in\CC(S)}
			{\max_{Q\in\QQ}\lambda_{Q,\sigma}(\alpha) \over 
				 \max_{Q\in\QQ}\lambda_{Q,\tau}(\alpha)} \nonumber\\
		&\le &
		\sup_{\alpha\in\CC(S)} 
    		{\lambda_{Q_\alpha,\sigma}(\alpha)\over  
		 \lambda_{Q_\alpha,\tau}(\alpha)} \nonumber\\
		&\le &
		\max_{Q\in\QQ} \sup_{\alpha\in\CC(S)} 
				{\lambda_{Q,\sigma}(\alpha)\over   
				 \lambda_{Q,\tau}(\alpha)}
\label{ratio le}
\end{eqnarray}

For the inequality in the other direction, 
consider first a non-annular component $Q\in\QQ$.
By definition (\ref{P length def}) of $\lambda_{Q,\sigma}$ the
supremum 
$\sup_{\alpha\in\CC(S)} 
				{\lambda_{Q,\sigma}(\alpha) / 
				 \lambda_{Q,\tau}(\alpha)}
$
may be obtained by letting $\alpha$ vary over 
$\CC(Q,\boundary Q)$.  By lemma \ref{closed curves give ratio}, we
may further  
restrict to $\alpha\in\CC_0(Q)$, at the expense of a bounded factor.
By lemma \ref{insulated surface},
$\lambda_{Q,\sigma}(\alpha)\simmult\lambda_{S,\sigma}(\alpha)$ 
for $\alpha\in\CC_0(Q)$ and any $\sigma\in Thin_\ep(S,\gamma)$.
Since $\CC_0(Q)\subset\CC(S)$, the supremal ratio over all of $\CC(S)$ can 
only be bigger, so in this case we have the desired inequality
\begin{equation}
\sup_{\alpha\in\CC(S)} {\lambda_{S,\sigma}(\alpha)\over
			\lambda_{S,\tau}(\alpha)}
\ge 
\sup_{\alpha\in\CC(S)} {\lambda_{Q,\sigma}(\alpha) \over
			\lambda_{Q,\tau}(\alpha)  }.
\label{ratio ge Q}
\end{equation}

%
%
%

Now consider an annulus $A\in\QQ$, and let $\alpha$ be a curve nearly 
maximizing the ratio
$\lambda_{A,\sigma}(\alpha)/\lambda_{A,\tau}(\alpha)$. By lemma
\ref{closed curve from annulus segment}, there exists
$\alpha'\in\CC(S)$ such that
\begin{equation}
\label{same twist}
|t_{A,\tau}(\alpha') - t_{A,\tau}(\alpha)|\le c_1,
\end{equation}
and such  that (applying theorem \ref{main estimate})
\begin{equation}
\label{annulus dominates}
\lambda_{S,\tau}(\alpha') \simmult \lambda_{A,\tau}(\alpha').
\end{equation}
Applying lemmas \ref{internal and external
twists} and \ref{compare twists}, we see that 
$t_{A,\tau}(\alpha') - t_{A,\sigma}(\alpha')$ and
$t_{A,\tau}(\alpha) - t_{A,\sigma}(\alpha)$  are both approximated, up
to bounded additive error, by the twist-parameter difference 
$s_j(\tau)-s_j(\sigma)$ associated to $A$. This together with inequality
(\ref{same twist}) implies that 
$|t_{A,\sigma}(\alpha') - t_{A,\sigma}(\alpha)|$ is uniformly
bounded as well.  Since $\lambda_{A,\sigma}$ and $\lambda_{A,\tau}$
are determined by the twisting numbers, we conclude that
$$
{\lambda_{A,\sigma}(\alpha') \over \lambda_{A,\tau}(\alpha')} \simmult
{\lambda_{A,\sigma}(\alpha) \over \lambda_{A,\tau}(\alpha)}.
$$
Now we may apply (\ref{annulus dominates}) together with the inequality
$\lambda_{A,\sigma}(\alpha') \lesim \lambda_{S,\sigma}(\alpha')$,
which is a consequence of theorem \ref{main estimate},
to conclude that 
\begin{equation}
{\lambda_{S,\sigma}(\alpha') \over \lambda_{S,\tau}(\alpha')} \gesim
{\lambda_{A,\sigma}(\alpha') \over \lambda_{A,\tau}(\alpha')}. 
\label{ratio ge A}
\end{equation}
Inequalities (\ref{ratio le}), (\ref{ratio ge Q}) and (\ref{ratio ge A})
together imply (\ref{Reduce to relative subsurface}).

\medskip

The length ratio on the left side of (\ref{Reduce to relative
subsurface}) is, by Kerckhoff's theorem \ref{Kerckhoff distance thm}, 
equal to $\exp 2d_{\TT(S)}(\sigma,\tau)$.
It remains to connect the length ratios on the right of  (\ref{Reduce
to relative subsurface}) to the distances 
$d_i \ \ (i=0,\ldots,k)$, so that we 
may complete the proof of 
theorem \ref{Product structure}.

Consider first $d_0$, the Teichm\"uller distance on $\TT(S_\gamma)$. 
Given a $Q\in\QQ$ which is a component of $S_\gamma$, with 
the conformal structure $\Pi_0(\sigma)$  or $\Pi_0(\tau)$ (restricted
from $S_\gamma$ to $Q$), we have by Kerckhoff's theorem:
$$\sup_{\alpha\in\CC_0(Q)}{\lambda_{Q,\Pi_0(\sigma)}(\alpha)\over 
			\lambda_{Q,\Pi_0(\tau)}(\alpha)} 
			= \exp 2d_{\TT(Q)}(\Pi_0(\sigma),\Pi_0(\tau)).
$$

We now claim that the extremal length ratio on the left can be
approximated by a ratio taken with respect to $\sigma$ and $\tau$. In
fact we claim that $(Q,\sigma)$ embeds $K$-quasi-conformally, with
uniform $K$, in $(Q,\Pi_0(\sigma))$. Recall the construction: we
measure the length and twist parameters for
$\gamma_{k+1},\ldots,\gamma_n$ for $\sigma$, and then
construct $\Pi_0(\sigma)$ with the same parameters, but with the
boundary curves $\gamma_1,\ldots,\gamma_k$ replaced by punctures. Thus
the embedding of $(Q,\sigma)$ in $(Q,\Pi_0(\sigma))$ can be done
separately on each pair of pants. After splitting each pair of pants
along seam curves into right-angled hexagons, the problem reduces to
the following: Let $H_1$ be a hyperbolic right-angled hexagon with lengths
$\ell_1,\ell_2,\ell_3$ on alternating sides, where $\ell_1\le \ep_1$,
and let $H_2$ be a hyperbolic right-angled hexagon with lengths
$0,\ell_2,\ell_3$ (really a right-angled pentagon with one ideal
vertex). There is a $K$-quasiconformal embedding (in fact 
bilipschitz) of $H_1$ into $H_2$ which is an isometry on sides $2$ and
$3$, where $K$ depends on $\ep_1$, but not on the $\ell_i$. 
This fact follows easily from the hexagon cosine law (see e.g.
\cite{beardon}). 
A similar computation applies when a pair of pants has more than one
short boundary which is to be replaced by punctures. The maps on the
pairs of pants fit together across the boundaries, because the marking
and twist parameters are the same for both surfaces.

It now follows (applying also lemma \ref{insulated surface}) that 
$\lambda_{Q,\sigma}(\alpha) \simmult \lambda_{Q,\Pi_0(\sigma)}(\alpha)$
uniformly for all $\sigma\in Thin_\ep(S,\gamma)$ and $\alpha\in\CC_0(Q)$. 
Thus
$$
\sup_{\alpha\in\CC_0(Q)}{\lambda_{Q,\Pi_0(\sigma)}(\alpha)\over
					\lambda_{Q,\Pi_0(\tau)}(\alpha)} 
\simmult
\sup_{\alpha\in\CC_0(Q)}{\lambda_{Q,\sigma}(\alpha)\over
						\lambda_{Q,\tau}(\alpha)}.
$$
Applying again lemma \ref{closed curves give ratio}, we conclude
$$
\exp 2d_{\TT(Q)}(\Pi_0(\sigma),\Pi_0(\tau)) \simmult
  \sup_{\alpha\in\CC(Q)}{\lambda_{Q,\sigma}(\alpha)\over
						\lambda_{Q,\tau}(\alpha)}.
$$
Since $d_0=d_{\TT(S_\gamma)}$ is the supremum of the distances in the 
Teichm\"uller spaces of the components, we have shown that the 
contribution to the right side of (\ref{Reduce to relative subsurface}) 
coming from the non-annular components is, up to bounded factor,  just
$\exp 2d_0$.  

We next consider an annular component $A_j$ in $\QQ$, and estimate its
extremal length contribution $\lambda_{A_j,\sigma}(\alpha)$, for
$\alpha\in\CC(S)$ and $\sigma\in Thin_\ep(S,\gamma)$, in terms of the
coordinates $\Pi_j(\sigma) = (s_j(\sigma),1/\ell_\sigma(\gamma_j))$. 
Rewrite $\Pi_j(\sigma) = (x_j(\sigma),y_j(\sigma))$, and
note first that (by (\ref{modulus for short geodesic}))
$y_j(\sigma)$ is closely approximated by
$m(A_j,\sigma)/\pi$.

Fix some arbitrary $\sigma_0\in Thin_\ep(S,\gamma)$ and define
$b_j(\alpha) = t_{\gamma_j,\sigma_0}(\alpha) - s_j(\sigma_0)$.
By lemma \ref{compare twists}, 
$t_{\gamma_j,\sigma}(\alpha)$ is estimated by  $b_j(\alpha) +
s_j(\sigma)$, up to bounded additive error, for all $\sigma$ and $\alpha$.
Applying also  lemma \ref{internal and external twists}, this
gives an estimate for $t_{A_j,\sigma}(\alpha)$. 
Now 
applying definition (\ref{A length def}) when $i(\alpha,\gamma_j) >
0$, we obtain this approximation for the length ratio contribution of
the annulus:  
\begin{equation}
 {\lambda_{A_j,\sigma}(\alpha) \over
                        \lambda_{A_j,\tau}(\alpha)} 
\simmult                        
 {y_j(\sigma) + (b_j(\alpha) +
				x_j(\sigma))^2/y_j(\sigma) \over 
                        y_j(\tau) + (b_j(\alpha) + 
                                       x_j(\tau))^2/y_j(\tau)},
\label{annulus length ratio}
\end{equation}
Note that the case where $i(\alpha,\gamma_j)=0$ and
$n(\alpha,\gamma_j)>0$ can be obtained by the limit of (\ref{annulus
length ratio}) as $b_j\to \infty$. 

By lemma \ref{closed curve from annulus segment},
we know that $t_{\gamma_j,\sigma_0}(\alpha)$, and therefore
$b_j(\alpha)$, can take on any value in $\R$ up to bounded difference. 
Since $y_j(\sigma),y_j(\tau)$ are at least $1/\ep$, we may once again
convert the additive error to a multiplicative one. Thus, taking a supremum
of (\ref{annulus length ratio}) over all $\alpha\in\CC(S)$, and
applying lemma \ref{Kerckhoff-Royden for torus}, we obtain
$$
\sup_{\alpha\in\CC(S)} {\lambda_{A_j,\sigma}(\alpha) \over
                        \lambda_{A_j,\tau}(\alpha)} 
\simmult
\exp 2d_{\Hyp_j}(\Pi_j(\sigma),\Pi_j(\tau)).
$$ 
It follows that the right side of (\ref{Reduce to relative subsurface})
is, up to bounded multiple, 
$$\max_{j=0,\ldots,k} \exp
2 d_j(\Pi_j(\sigma),\Pi_j(\tau)).
$$
Taking logarithms, we obtain the desired statement of theorem
\ref{Product structure}.
\end{proof}

\section{Geometric remarks}
\label{remarks}
Having shown that regions of the form $Thin_\ep(S,\gamma)$ in $\TT(S)$
are approximated up to bounded  additive distortion
by products with the sup metric, we
discuss briefly here the implications to ``non-hyperbolic'' behavior
of the Teichm\"uller space. Indeed, a product with the sup metric is
not just non-negatively curved, but exhibits some definite
positive-curvature behavior, as we shall see.


Let $X$ be a metric space, and let $x,y,z\in X$ denote three points. 
We say that $z$ is
``$\delta$-between'' $x$ and $y$ if 
$$ d(x,z) + d(z,y)  - d(x,y) < \delta$$
for $\delta>0$. In particular if $X$ is a length space and $[xy]$
denotes a shortest path (geodesic) between $x$ and $y$, then when $z$
is $\delta$-between $x$ and $y$, the path
$[xz]\union[zy]$ is a
quasigeodesic with only additive distortion $\delta$.

In a negatively curved length space in the sense of Gromov and Cannon,
if $z$ is $\delta$-between $x$ and $y$ 
then it is contained in an $R(\delta)$-neighborhood of $[xy]$, where
$R(\delta)$ is independent of $x$ and $y$. This is
called stability of quasi-geodesics. 

Let us define an ``instability function'' $s$ for $X$ as follows: For
$\delta,L\ge 0$ define $s(\delta,L)$ by
$$s(\delta,L) = \sup\{d(z,[xy]): \text{$z$ is  $\delta$-between $x$
and $y$ and  $d(x,y)\le L$}\}.$$
In order to consider only the large-scale behavior, let us denote by
$\bar s(\delta)$ the asymptotic growth rate of $s(\delta,L)$; that is,
$$\bar s(\delta) =  \limsup_{L\to\infty} {\log s(\delta,L)\over \log L}.$$
It is easy to see that $\bar s(\delta)\in[0,1]$, and $\bar s$ is a
non-decreasing function of $\delta$. 

In a negatively curved space, $s(\delta,L)$ is independent of $L$ for
large $L$. Thus, $\bar s \equiv 0$.
In Euclidean space it is not hard to see that $s(\delta,L) =
\sqrt{2L\delta+\delta^2}/2$. In particular $\bar s(\delta) = 1/2$ for
$\delta>0$. 
In a space which contains large positively curved regions such as
hemispheres,  $s$ could be made to grow roughly proportionally to $L$, for
fixed $\delta$. Thus $\bar s(\delta)$ can be made identically $1$ in
such a space.   

Let $X$ denote the product $\R\times\R$ with the sup metric. Then even 
$s(0,L)$ in this space is at least $L/2$, as one can see by
considering $x=(0,0), y=(L,0)$ and $z = (L/2,L/2)$. 
Thus the growth rate $\bar s \equiv 1$, as in the case of positive curvature. 

Let $X,Y$ be two length spaces with instability functions $s_X$,
$s_Y$.  If $f:X\to Y$ is a map such that $|d_X - d_Y\circ f| \le c$,
it is easy to show that 
$$s_X(\delta,L) \le 3c + 4s_Y(\delta+3c,L+c),$$
and vice versa if $f$ is invertible. 
Applying this to the homeomorphism $\Pi$ of theorem \ref{Product
structure}, we may conclude that $\bar s_{\TT(S)}(\delta) =1$ for all
$\delta>0$. This is the
sense in which we say that $\TT(S)$ is positively curved in the large.

\ifx\undefined\bysame
\newcommand{\bysame}{\leavevmode\hbox to3em{\hrulefill}\,}
\fi

\end{document}